\input amstex


\def\b1{\text{\bf 1}}

\def\CC{{\Cal C}}
\def\CD{{\Cal D}}
\def\CF{{\Cal F}}
\def\CE{{\Cal E}}
\def\CH{{\Cal H}}

\def\CJ{{\Cal J}}
\def\CG{{\Cal G}}
\def\CL{{\Cal L}}
\def\CM{{\Cal M}}

\def\CO{{\Cal O}}
\def\CP{{\Cal P}}

\def\CS{{\Cal S}}

\def\CV{{\Cal V}}

\def\Hom{\text{Hom}}
\def\Sym{\text{Sym}}

\def\#{\,\check{}}

\def\fh{{\frak h}}

\def\fg{{\frak g}}
\def\fF{{\frak F}}

\def\ft{{\frak t}}
\def\fm{{\frak m}}
\def\fn{{\frak n}}

\def\fs{{\frak s}}

\def\Coker{\text{Coker}}
\def\Ker{\text{Ker}}

\def\Res{\text{Res}}
\def\Spec{\text{Spec}}

\def\limleft{\mathop{\vtop{\ialign{##\crcr
  \hfil\rm lim\hfil\crcr
  \noalign{\nointerlineskip}\leftarrowfill\crcr
  \noalign{\nointerlineskip}\crcr}}}}
\def\limright{\mathop{\vtop{\ialign{##\crcr
  \hfil\rm lim\hfil\crcr
  \noalign{\nointerlineskip}\rightarrowfill\crcr
  \noalign{\nointerlineskip}\crcr}}}}

\def\hra{\hookrightarrow}
\def\iso{\buildrel\sim\over\rightarrow}

\parskip=6pt

\documentstyle{amsppt}
\document
\magnification=1100
\NoBlackBoxes

\bigskip

\centerline {\bf Langlands parameters for
  Heisenberg modules.} 

\bigskip

\centerline {A.~Beilinson}

\bigskip

\centerline{\bf  Introduction} 

Below we
define a ``spectral decomposition" of the category of
representations of the Heisenberg Lie algebra; the
spectral parameters are moduli of de Rham local
systems. This is a very particular case (the group $G$ is
a torus) of the following  general conjecture  on 
Langlands parameters in the de Rham-Kac-Moody
(split) setting. The key device
here is chiral  Hecke  algebra (we assume its definition
is known to the reader; otherwise he should either consider
it as a black box or simply skip the introduction). When
$G$ is a torus the chiral Hecke algebra  is the same as
a usual lattice Heisenberg vertex algebra.

General idea: Representations of a p-adic group have a
canonical ``spectral decomposition" over the spectrum of
 Bernstein's center. According to the local Langlands
conjecture these parameters can be identified roughly with
the moduli  of Galois representations with values in
the Langlands dual group. A possible imitation of this
picture in the de Rham setting: replace the p-adic group
by a Kac-Moody Lie algebra\footnote{In truth
representations of the p-adic group correspond {\it not}
to individual Kac-Moody modules but rather
to categories of those invariant with respect to
the $G(F)$-action, see Remark (i) below.} and the Galois
representation by a de Rham local system. The fun is that
the usual Bernstein center is trivial here, which
corresponds to the fact that the moduli of  de Rham
local systems have no non-constant global functions. In
the de Rham setting however, contrary to the $p$-adic
situation, there is a direct way  to relate the Galois
side of the picture with ``automorphic" one provided by
the Hecke chiral algebras.

One proceeds as follows. Below $F\simeq k((t))$ is
a local $k$-field,
$O\subset F, O\iso k[[t]]$ the local ring (here $k$ is our
base field of char 0). Let $G$
be a reductive group over $k$, $c$ an integral level
for $G$ which is less than critical. Let
$\CH^c$ be the chiral Hecke algebra of level $c$. The
Langlands dual group $G^\vee$  acts on $\CH^c$,
so you can twist $\CH^c$ by any $ G^\vee$-local system
$\phi$ (in the de Rham sense) on Spec$F$ getting a chiral
algebra
$\CH_\phi^c$ on Spec$F$.   When
$\phi$ varies we get a family of chiral algebras
$\CH_{\CL\CS}^c$ parametrized by the moduli stack
$\CL\CS = \CL\CS_{ G^\vee} $ of  $ G^\vee$-local
systems on Spec$F$. Consider the
 category of $\CH_{\CL\CS}^c$-modules.\footnote{supported
at the puncture of Spec$F$.} Every such animal $V$ is an
$\CO$-module on $\CL\CS$ equipped with an
extra structure. In particular, it carries an action of the
``constant" (with respect to $\phi$) vertex algebra
 of $ G^\vee$-invariants in $\CH^c$,
which is the enveloping algebra of the Kac-Moody Lie
algebra $\fg (F)^c$ of level
$c$. Thus
$\Gamma (\CL\CS ,V)$ is a $\fg (F)^c$-module. 

{\bf Conjecture.} The functor $\Gamma (\CL\CS,\cdot ):
\{ \CH_{\CL\CS}^c$-modules$\} \to \{ \fg
(F)^c$-modules$\}$  is an equivalence of
categories.\footnote{We consider $\fg (F)^c$ as a
topological Lie algebra, so
$\fg (F)^c$-modules are assumed to have the property that
every vector is killed by $t^n \fg (O) \subset \fg (F)^c$
for sufficiently large $n$. The central
element $1\in k\subset
\fg (F)^c$  acts as identity.} 

Roughly speaking, this means that
$\CH^c_\phi$-modules equipped with an action of
Aut$\phi$ are the same as certain special $\fg
(F)^c$-modules; 
$\fg (F)^c$-modules coming from $\CH^c_\phi$'s with
different
$\phi$'s are different; and every $\fg (F)^c$-module can
be presented as an ``integral" of modules coming from
$\CH^c_\phi$-modules. Of course, a precise formulation
of the above conjecture requires an explanation of what
are $\CO$-modules on $\CL\CS$ (which is an ind-algebraic
stack in a broad sense), etc.

{\it Remarks.} (i) The chiral algebra $\CH^c$ on Spec$\,
F$ carries a natural $G(F)$ action commuting with the
$G^\vee$-action, so $G(F)$ acts on $\CH_{\CL\CS}^c$ as
well. Thus $G(F)$ acts on the category of
$\CH_{\CL\CS}^c$-modules; it acts also on the category of
$\fg (F)^c$-modules (via the adjoint action of $G(F)$ on
$\fg (F)^c$). The functor $\Gamma (\CL\CS,\cdot )$
commutes with the $G(F)$-action.

(ii) Assuming the above conjecture, one  defines
the (Galois) support of a $\fg (F)^c$-module as the support
of the corresponding $\CH_{\CL\CS}^c$-module. It would be
very interesting to describe explicitly the categories of
$\fg (F)^c$-modules with  support in a given subspace of $
\CL\CS$.  For example, one can hope that the
category of $\fg (F)^c$-modules supported on the substack
of regular singular connections with nilpotent residue  is
equivalent to the product of several copies of the category
of $\CD$-modules on the affine flags space $\CF_G$ of
$G(F)$.\footnote{ $\CF_G$ is a formally smooth ind-proper
ind-scheme.} The copies are labeled by the
orbits of the
$c$-affine action of the affine Weyl group on
lattice of weights. As in the finite-dimensional
situation, the functor that coorrresponds to an
orbit assigns to a
$\CD$-module $M$ the vector space of global sections
$\Gamma_\CL (M):=\Gamma (\CF_G ,M\otimes\CL )$ where $\CL$
is the (positive) line bundle on $\CF_G$ that corresponds
to the (only) negative weight from the orbit; this is a
$\fg (F)^c$-module in an obvious way. It is known that the
functor $\Gamma_\CL$ is exact and fully faithful; a
construction due to Gaitsgory provides a canonical lifting
of $\Gamma_\CL$ to a functor with values in
$\CH^c_{\CL\CS}$-modules. More generally, the category of
$\fg (F)^c$-modules supported on the substack of 
connections with regular singularities with fixed
``semi-simple part of the monodromy"  should be equivalent
to the product of several copies of the category of
appropriately twisted
$\CD$-modules on
$\CF_G$.  It would be very interesting to guess  what
happens in the case of irregular singularities.

(iii) Denote by $T$ the Cartan group of $G$, so
$T^\vee \subset G^\vee$ is the Cartan torus of
$G^\vee$. Let
$V$ be a
$\fg (F)^c$-module,
$\phi$ a
$G^\vee$-local system which comes from a
$T^\vee$-local system $\phi_{T^\vee}$. One can hope
that $\phi$ is in the support of $V$ if and only if
$\phi_{T^\vee}$ belongs to the support of one of the
Heisenberg modules $H^{\infty /2 +\cdot}(\fn (F),V)$ where
$\fn$ is the Lie algebra of a maximal nilpotent subgroup
$N\subset G$.

\medskip

In this note we check the toy case of the above conjecture
when
$G=T$ is a torus, hence the Hecke chiral algebra is the
lattice Heisenberg chiral algebra. The result is
essentially straightforward (the key point is that
 $\CL\CS_{T^\vee}$ is covered by $\Spec (\Sym
\ft (O))$ = the spectrum of the algebra of
``annihilation operators" in Heisenberg Lie algebra
modules); the pages below are mostly a review of the known
definitions and constructions.
\S 1 considers a class of Heisenberg groups we play with.
In \S 2  we transplant the picture to a curve.  \S 3
deals with lattice vertex algebras. In case of a
non-degenerate level we establish a natural bijection
between symmetric Heisenberg groups and symmetric lattice
vertex algebras (the construction is essentially the same
as in [K] 5.4, 5.5 though the presentation may look
different) and identify the corresponding categories of
representations  (the latter subject is essentially
contained in [D], the proof below is similar to
Dong's\footnote{D.~Gaitsgory found recently another, more
natural, proof.}). We also show that the
twist by a local system of a lattice chiral algebra
 does not affect (locally) the chiral algebra
structure but changes the Heisenberg Lie algebra
embedding. In \S 4 the structure of the moduli space of
local systems is spelled out and the promised fact on
Langlands parameters is stated and proved.

I am grateful to V.~Drinfeld, D.~Gaitsgory and D.~Kazhdan
for their interest and stimulating discussions.

\bigskip

\centerline{\bf \S1  Symmetric Heisenberg extensions} 

\medskip 

{\bf 1.1} Suppose $A$, $C$ are abelian groups,
$H$ a central $C$-extension of $A$. Denote by
$\{ \, ,\, \}$ the corresponding commutator pairing
$A\times A\to C$. For $a\in A$ we denote its preimage in
$H$ by $H_a$; this is a $C$-torsor. We write elements of
$H_a$  as $\tilde{a}$.

A {\it symmetric structure} on $H$ is an
automorphism
$\sigma$ of $H$ such that $\sigma^2
=$id$_{H}$,
$\sigma|_C =$id$_C$, and $\sigma$mod$C$ is the involution
$a\mapsto a^{-1}$ of $A$. We call  $(H,\sigma )$ a {\it
symmetric (central) $C$-extension} of $A$. 

{\it Example.} Let $s :A\to H$ be a set-theoretic section
such that the corresponding cocycle $a_1 ,a_2 \mapsto
\gamma_s (a_1 ,a_2 ):= s(a_1 )s(a_2 )s(a_1 a_2 )^{-1}$ has
property $\gamma_s (a_1^{-1}, a_2^{-1})=\gamma_s (a_1 ,a_2
)$. Then $s$ defines a symmetric structure such that
$\sigma s(a)=s(a^{-1})$. In other words, the extension
defined by a two cocycle with the above property is
symmetric.

 All central $C$-extensions of $A$ form a groupoid. The
Baer
product defines on it the structure of a Picard groupoid.
Same is true for the groupoid of symmetric extensions. The
 automorphism group of a cental extension  equals
Hom$(A,C)$;\footnote{We identify $f:A\to C$ with the
automorphism $\tilde{a}\mapsto \tilde{a}f(a)$.}  the
one of a symmetric extension is Hom$(A,C_2
)\subset$ Hom$(A,C)$.\footnote{ Here
$C_2
\subset C$ is the subgroup of elements of order 2.}

\medskip

{\bf 1.2} Denote by
$H^{(2)}$ the Baer product of two copies of $H$. We
have a canonical homomorphism $H\to H^{(2)}$, $h\mapsto
h^{(2)}$, which lifts id$_A$ and such that $c^{(2)}=c^2$
for $c\in C$. A symmetric structure
$\sigma$ on
$H$ yields a set-theoretic section
$s_\sigma : A\to H^{(2)}$. Namely, for $a\in A$  our
$s_\sigma (a)$ is the trivialization of the $C$-torsor
$H^{(2)}_a$ defined by the identification $H^{(2)}_a :=
H_a \cdot H_a \iso C$, $h\cdot h' \mapsto h\sigma
(h')=\sigma (h) h' \in C\subset H$.  

{\bf Lemma.} The map $\sigma \mapsto s_\sigma$ is a
bijection between the set of symmetric structures on $H$
and the set of set-theoretic sections $s: A\to H^{(2)}$
which satisfy an equation $s(a)s(b)=s(ab)\{ a,b\}$.  

{\it Proof.} One has $s_\sigma
(a)=[\tilde{a}\sigma
(\tilde{a})]^{-1}\tilde{a}^{(2)}$ where $\tilde{a}\in H$
is a lifting of $a$ and $[\tilde{a}\sigma
(\tilde{a})]\in C$. Therefore
$s_\sigma (a) s_\sigma (b)=[\tilde{a}\sigma
(\tilde{a})]^{-1}\tilde{a}^{(2)}[\tilde{b}\sigma
(\tilde{b})]^{-1}\tilde{b}^{(2)}= [\tilde{a}\sigma
(\tilde{a})\tilde{b}\sigma
(\tilde{b})]^{-1}
(\tilde{a}\tilde{b})^{(2)}=\{
a,b\}[(\tilde{a}\tilde{b})\sigma (\tilde{a}\tilde{b})]^{-1}
(\tilde{a}\tilde{b})^{(2)}=\{ a,b\} s_\sigma (ab),$ so
$s_\sigma$ satisfies our equation. Conversely, for a
given $s$ we define  $\sigma_s :H\to H$ as $\sigma_s
(\tilde{a}):= [\tilde{a}^{(2)}s(a)^{-1}] \tilde{a}^{-1}$
where $[\tilde{a}^{(2)}s(a)^{-1}]\in C$. We leave it to
the reader to check that
$\sigma
$ is a symmetric structure on $H$.
\hfill$\square$

In other words, for a fixed skew-symmetric pairing $\{ \,
,\,\} :C\times C\to A$ a symmetric
$C$-extension with the commutant pairing equal to $\{ \,
,\,\}$ is the same as a Baer square root of the
 $C$-extension of $A$ defined by the 2-cocycle 
$\{ \,
,\,\}$.

\medskip

{\bf 1.3} We see that the Picard groupoid of 
 commutative symmetric $C$-extensions of $A$  is
canonically equivalent to the Picard groupoid $\CE xt (A,
C_2 )$ of commutative $C_2$-extensions of
$A$. If
$A$ is a free commutative group then this is the category
of Hom$(A,C_2 )$-torsors.

 For a fixed pairing $\{ \, ,\,\} :A\times A \to C$  the
groupoid of symmetric
$C$-extensions of $A$ having $\{
\, ,\, \}$ as the commutator pairing carries the Baer sum
action of the Picard groupoid of commutative symmetric
extensions, i.e., that of $\CE xt (A,
C_2 )$. If our groupoid  is non-empty then it is an 
$\CE xt (A, C_2 )$-torsor.

\medskip

 The above considerations make sense if $H$ is a
central extension of an abelian group $A$ by any Picard
category; in particular, we can play with superextensions
of
$A$.\footnote{See [BBE] \S2 for the terminology.} The
above lemma remains true in this more general setting.

\medskip

 {\bf 1.4} All schemes  below are over the base field $k$
of char 0. We will consider symmetric $\Bbb G_m$-extensions
or superextensions of certain group ind-schemes. 

{\it Remark.} Let $A$ be a group ind-scheme equipped
with a pairity homomorphism $A\to \Bbb Z/2$ and a (super)
skew-symmetric  pairing
$\{ \, ,\, \} :A\times A \to \Bbb G_m$. As follows from
1.2, symmetric superextensions of $A$ having $\{ \, ,\,
\}$ as the commutator pairing are in bijective
correspondence with those for the
reduced ind-scheme $A_{red}$.

Let $F\simeq k((t))$ be a local $k$-field,
$O\simeq k[[t]]$ its local ring, $\fm_x \simeq tk[[t]]$
the maximal ideal, $x\in$ Spec$(O)$ is the closed point. 

Let
$T=\Gamma
\otimes\Bbb G_m$ be a torus, $\ft =\Gamma \otimes k$ its
Lie algebra. We fix a  level $c$ which is an
integral symmetric bilinear form
on $\Gamma$.  We have the
commutative group ind-scheme $T(F)=\Gamma \otimes
F^\times$. For $\gamma\in\Gamma$, $f\in F^\times$ we set
$f^\gamma :=\gamma \otimes f\in T(F)$. The Lie algebra of
$T(F)$ equals $\ft (F)=\Gamma \otimes F$; for
$\gamma\in\Gamma$,
$\varphi \in F$ we set $\varphi\cdot \gamma :=
\gamma\otimes
\varphi \in \ft (F)$.  The valuation  $v : F^\times
\to
\Bbb Z$ yields a projection $T(F)\twoheadrightarrow
\Gamma$;  for
$\gamma \in
\Gamma$ the preimage of $\gamma$ is
$T(F)^\gamma \subset T(F)$. Notice that 
$T(F)^0_{red}$ is the group scheme $T(O)$.

{\bf  Definition.} A {\it  Heisenberg extension of level
$c$} is a central $\Bbb G_m$-superextension $
T(F)\,\tilde{}$ of $T(K)$  equipped with a splitting
$i: T(O)\to T(F)\,\tilde{}$ which satisfies the following
two properties:

(i) The $
T(F)\,\tilde{}\,$-pairity of
$a\in T(F)_\gamma$ equals $c(\gamma ,\gamma )\,$mod$2$. 

(ii) The
 $
T(F)\,\tilde{}\,$-commutator pairing
$ T(F)\times T(F)\to \Bbb G_m$ for $
T(F)\,\tilde{}$ is $$f_1^{\gamma_1},
f_2^{\gamma_2} \mapsto \{ f_1^{\gamma_1}
,f_2^{\gamma_2} \}^c := \{ f_1 ,f_2
\}^{-c(\gamma_1 ,\gamma_2 )} \tag 1.4.1$$ where $ \{ f_1
,f_2
\}\in \Bbb G_m$ is the Contou-Carr\`ere
symbol.\footnote{See [CC] or [BBE]
\S 3.}

We say that $
T(F)\,\tilde{}$ is  {\it symmetric} it is equipped with a
symmetric structure $\sigma$ which fixes the splitting $i$.

{\bf 1.5 Remarks.} (i) The Baer product of
(symmetric) Heisenberg extensions of levels $c, c'$ is a
(symmetric) Heisenberg extension of level $c+c'$.

(ii) The Lie algebra $\ft (F)\,\tilde{}$ of a
Heisenberg extension of level $c$ is a central
$k$-extension of the commutative Lie algebra $\ft (F)$.
The corresponding  commutator pairing $\ft
(F)\times\ft (F)
\to k$ is
 $\varphi_1 \cdot \gamma_1 ,\varphi_2 \cdot\gamma_2
\mapsto c(\gamma_1 ,\gamma_2 )$Res$(\varphi_2 d\varphi_1
)$. A symmetric structure defines a $k$-linear continuous
splitting $\tilde{i} : \ft (F) \to \ft (F)\,\tilde{}$ (the
only splitting fixed by $\sigma$) which equals $i$ on $\ft
(O)$. Thus $\ft (F)\,\tilde{}$ is uniquely defined: as a
vector space it equals $\ft (F)\oplus k$ with bracket
defined by the above formula. We call $\ft (F)\,\tilde{}$
{\it the} symmetric Heisenberg Lie algebra of level $c$.
The adjoint action of the group $T(F)$ on $\ft
(F)\,\tilde{}$ is Ad$_{f^\gamma}(\varphi \cdot \gamma' )=
\varphi \cdot \gamma' + c(\gamma ,\gamma' )$Res$ (\varphi
d\log f)$.

(iii) Let $T^\vee =\Gamma^\vee \otimes\Bbb G_m$ be the
dual torus. The group ind-scheme $T^\vee (F)$ acts on
$T(F)\,\tilde{}$ according to formula $g^{\check{\gamma}}
(\tilde{f}^\gamma )=\{ f,g\} ^{(\check{\gamma},\gamma
)}\tilde{f}^\gamma$ where $f,g \in F^\times$, $\gamma
\in\Gamma$, $\check{\gamma}\in \Gamma^\vee$, and
$\tilde{f}^\gamma \in T(F)\,\tilde{}$ is a lifting of
$f^\gamma$. One can view $c$ as a homomorphism $T\to
T^\vee$. Then for any $f^\gamma \in T(F)$ the adjoint
action of $f^\gamma$ on $T(F)\,\tilde{}$ coincides with
the action of $c(f^\gamma )\in T^\vee (F)$.

(iv) As follows from 1.2 for a symmetric Heisenberg
extension
$T (F)\,\tilde{}$ the action of the group ind-scheme
Aut$(F)$ on $T(F)$ lifts uniquely to an action on $T
(F)\,\tilde{}$ of the two-sheeted covering Aut$^{1/2}(F)$
which preserves $\sigma$ and $i$. In particular, $T
(F)\,\tilde{}$ carries a canonical action of the Lie
algebra $\Theta (F):=$ Der$(F)$ of vector fields on $F$.
Here is an explicit formula:

{\bf 1.6 Lemma.} For $f\in F^\times$, $\gamma \in \Gamma$,
$\theta \in \Theta (F)$, and a
lifting $\tilde{f}^\gamma  \in T(F)\,\tilde{}$  one has
$$\theta (\tilde{f}^\gamma )= [\tilde{i}(\frac{\theta
(f)}{f}
\cdot \gamma )  +\frac{c(\gamma
,\gamma )}{2}\Res (\theta (f) \frac{df}{f^2})
]\tilde{f}^\gamma =\tilde{f}^\gamma
[\tilde{i}(\frac{\theta (f)}{f}
\cdot \gamma )-\frac{c(\gamma
,\gamma )}{2}\Res (\theta (f) \frac{df}{f^2}) ].$$

{\it Proof.} The action of Aut$^{1/2} (F)$ is compatible
with the projection $T(F)\,\tilde{} \to \Bbb G_m$,
$\tilde{f}^\gamma \mapsto \tilde{f}^\gamma \sigma
(\tilde{f}^\gamma )$, see 1.2. Therefore 
$(1+\epsilon \theta )\tilde{f}^\gamma =(1+\epsilon 
\tilde{i}(\frac{\theta (f)}{2f} \cdot
\gamma ))\tilde{f}^\gamma (1+\epsilon
\tilde{i}(\frac{\theta (f)}{2f} \cdot
\gamma ))$. Now use (1.4.1).  \hfill$\square$

Denote by $\CH s^c$ the groupoid of symmetric Heisenberg
extensions of level $c$. By 1.2 and
1.5(i) 
$\CH s^0$ is a Picard groupoid canonically equivalent to
that of Hom$(\Gamma ,\mu_2
)=\Gamma^\vee \otimes\mu_2$-torsors.\footnote{Here
$\Gamma^\vee$ is the dual lattice.} Thus for any $c$
the groupoid $\CH s^c$ is
a $\Gamma^\vee \otimes\mu_2$-gerbe (the fact that $\CH
s^c$ is non-empty follows from, say, Remark in 1.7 below). 

 {\bf 1.7} Here is anoher description of $\CH
s^c$. Let $T(F)\,\tilde{}$ be a symmetric Heisenberg
extension. For $\gamma \in \Gamma$ let $\lambda^\gamma$ be
 a superline defined as follows. Consider
$T (F)\tilde{}$ as a superline $\lambda$ over
$T(F)$ equivariant with respect to the left and right
$T(F)\,\tilde{}$-translations.  Restrict $\lambda$ to
$T(F)_{red}^{-\gamma}$ and consider the right
translation action of $T(O)\buildrel{i}\over\hra
T(F)\,\tilde{}$. Our $\lambda^\gamma$ is the superline of
$T(O)$-equivariant sections.\footnote{We use notation
$\lambda^\gamma$ instead of 
$\lambda^{-\gamma}$ since $T(O)$ acts on it by the
character $c(\gamma )$, see 3.2, 3.3.}

Denote by
$\omega_x$ the cotangent $k$-line $\fm_x /\fm_x^2$.

{\bf  Lemma.} (i) For every $\gamma_1 ,\gamma_2 \in
\Gamma$ there is a canonical isomorphism  
$\mu  =\mu_{\gamma_1 \gamma_2}:\lambda^{\gamma_1}
\otimes \lambda^{\gamma_2}\iso \lambda^{\gamma_1
+\gamma_2}\otimes
\omega_x^{c(\gamma_1 ,\gamma_2 )}$, $\ell_1 \otimes\ell_2
\mapsto
\mu (\ell_1 ,\ell_2 )$, and for every $\gamma \in\Gamma$ a
one
$\sigma =\sigma_\gamma :
\lambda^\gamma \iso \lambda^{-\gamma}$
such that

(a) $\mu$ is associative: for every $\gamma_1 ,\gamma_2
,\gamma_3 \in \Gamma$ and
$\ell_i \in \lambda^{\gamma_i}$ one has 
$\mu (\mu (\ell_1 ,\ell_2 ),\ell_3 )=
\mu (\ell_1 ,\mu (\ell_2 ,\ell_3 ))\in
\lambda^{\gamma_1 +\gamma_2 +\gamma_3}\otimes
\omega_x^{c(\gamma_1 ,\gamma_2 )+c(\gamma_2 ,\gamma_3
)+c(\gamma_1 ,\gamma_3)}$.

(b) One has $\mu_{\gamma_1
\gamma_2}=(-1)^{c(\gamma_1 ,\gamma_2 )}\mu_{\gamma_2
\gamma_1}\,\fs$ where $\fs :\lambda^{\gamma_1}
\otimes\lambda^{\gamma_2}\iso
\lambda^{\gamma_2}\otimes\lambda^{\gamma_1}$ is the
commutativity constraint.

(c) $\sigma^2 =$ id, i.e.,
$\sigma_{\gamma}=\sigma_{-\gamma }^{-1}$, and
$\mu$ commutes with $\sigma$.

(ii) The above construction identifies $\CH s^c$ with the
groupoid of triples $(\{ \lambda^\gamma \} ,\mu ,\sigma )$ 
where
$\lambda^\gamma$, $\gamma
\in\Gamma$, is a collection of superlines such that each
$\lambda^\gamma$ has pairity
$c(\gamma ,\gamma )$ mod2, and $\mu$, $\sigma$ are data of
isomorphisms as in (i) above.

{\it Remark} (cf.~[K] 5.5). Choose a square root
$\lambda_x$ of
$\omega_x$ and consider it as an odd superline. For
a datum as in Lemma set $\lambda^{'\gamma }:=
\lambda^\gamma \otimes \lambda_x^{\otimes c(\gamma
,\gamma )}$. Then $\mu$ and $\sigma$ make
$\{\lambda^{'\gamma}\}$  a symmetric $\Bbb G_m$-extension
of
$\Gamma$ with the commutator pairing $\gamma_1 ,\gamma_2
\mapsto (-1)^{c(\gamma_1 ,\gamma_2 )+c(\gamma_1
,\gamma_1 )c(\gamma_2 ,\gamma_2 )}$. Thus a choice of
$\lambda_x$ identifies
$\CH s^c$ with the category of such symmetric
$\Bbb G_m$-extensions of $\Gamma$.

{\it Proof of Lemma.} 
(i) The isomorphisms $\sigma$ come from the symmetric
structure on $T(F)\,\tilde{}$. In order
to define $\mu$ let us choose a parameter $t$ in $F$ (so
$t\in F$, $v (t)=1$). Then $\lambda^\gamma \iso
\lambda_{t^{-\gamma}}$, $\ell^\gamma \mapsto
\ell^{\gamma}_t := \ell^\gamma (t^{-\gamma} )$, and we 
define $\mu $ by formula $(\mu
(\ell^{\gamma_1} ,\ell^{\gamma_2}) (dt_0)^{
-c(\gamma_1 ,\gamma_2 )})_t := \ell^{\gamma_1}_{t}
\cdot\ell^{\gamma_2}_{ t}$ where $\cdot$ is the
product map in $T(F)\,\tilde{}$. It remains to check that
$\mu$ so defined does not depend on the choice of $t$.
Take another parameter
$at$,  $a\in O^\times$. Then $\ell^{\gamma}_{
at}=\ell^{\gamma}_{ t}\cdot i(a^{-\gamma} )$ hence
$\ell^{\gamma_1}_{at}\cdot \ell^{\gamma_2}_{ at}=
(\ell^{\gamma_1}_{t}\cdot\ell^{\gamma_2}_{t})\cdot
i(a^{-\gamma_1  -\gamma_2})\cdot \{ a ,t
\}^{-c(\gamma_1 ,\gamma_2
)}=(\ell^{\gamma_1}_{t}\cdot\ell^{\gamma_2}_{t})\cdot
i(a^{-\gamma_1 -\gamma_2})\cdot a_0^{-c(\gamma_1 ,\gamma_2
)}.$ Since
$d(at)_0 = a_0 dt_0$ this
implies the  independence. 

Properties (a)--(c) are immediate.

(ii) Assume we have a datum $(\{ \lambda^\gamma \}, \mu,
\sigma )$. By Remark in 1.4 a
symmetric Heisenberg extension of $T(F)$ amounts to a
similar extension of $T(F)_{red}$, so it suffices to
define $T(F)\tilde{}_{red}$. To do this we track back the
construction of (i). 
As a mere $T(O)$-equivariant superline
$T(F)\tilde{}_{red}$ is constant over each
$T(F)_{red}^\gamma$ with fiber $\lambda^{-\gamma}$. We also
have its trivialization $i$ over $T(O)$ and the symmetry
$\sigma$. To recover the product let us fix  a  parameter
$t$ in $F$. Then $\mu$ recovers the multiplication law over
the subgroup $t^\Gamma
\subset T(F)_{red}$ according to formula
$\ell^1_{t^{\gamma_1}}
\cdot
\ell^2_{t^{\gamma_2}}:= (\mu (\ell^1 ,
\ell^2 ) (dt_0 )^{-c(\gamma_1 ,\gamma_2
)})_{t^{\gamma_1 +\gamma_2}}$. It extends in a unique way
to the product over the whole
$T(F)_{red}$ if we demand that the $T(O)$-equivariant
structure coincides with the right multiplication via $i$
and the commutator pairing equals (1.4.1). The
independence of $t$ follows from a straightforward 
computation as in (i).
\hfill$\square$

{\bf 1.8 Example} (will not be used below). Assume that
$\Gamma =\Bbb Z$, so $T(F)=F^\times$, and $c=1$, i.e.,
$c(\gamma_1 ,\gamma_2 )=\gamma_1 \gamma_2$. Choose a
square root $\CL $ of the $F$-line $\omega (F)$. The group
$F^\times$ acts on the Tate $k$-vector space $\CL$ by
homotheties, so we have the Heisenberg group
$F^{\times\flat}_{(\CL )}$ defined as the pull-back of
the Tate superextension of the group GL$(\CL )$ of all
continuous $k$-automorphisms of $\CL$ (see e.g.~[BBE]
3.7). Now $F^{\times\flat}_{(\CL )}$ caries a canonical
symmetric structure. Indeed, the pairing
$\CL\times\CL \to k$, $\ell ,\ell' \mapsto
\text{Res}(\ell\ell' ),$ makes $\CL$ a self-dual Tate
$k$-vector space, and the symmetric structure $\sigma$
comes from the canonical isomorphism GL$^\flat (\CL
)\iso$ GL$^\flat (\CL )$ which lifts the automorphism
$g\mapsto {}^t g^{-1}$ (see [BBE] (2.16.1)). Notice that
the functor
$\CL\mapsto  F^{\times\flat}_{(\CL )}$ is a morphism,
hence an equivalence, of the $\mu_2$-gerbes. 

{\bf 1.9} We use notation of 1.5(iii). Set
$$Z^c :=
\Ker (c : T\to T^\vee 
)=\Hom (\Coker (c:\Gamma \to \Gamma^\vee ),\Bbb G_m ).
\tag 1.9.1$$ We consider $Z^c$ as a subgroup of $T\subset
T(O)$. The level
$c$ is non-degenerate if and only if $Z^c$ is finite.

Denote by $T(F)\,\tilde{}$mod the category of
$T(F)\,\tilde{}$-modules, i.e., $k$-vector (super)spaces
equipped with an action of $T(F)\,\tilde{}$ (the subgroup
$\Bbb G_m \subset T(F)\,\tilde{}$ is assumed to act by
standard homotheties) and by
$Z^c$mod the category of
$Z^c$-modules.

Assume  that $c$ is non-degenerate.
 
{\bf Lemma.} 
  Every set-theoretic section $s: \Gamma^\vee
/c(\Gamma )\to \Gamma^\vee$ defines an equivalence of
categories $e_s : T(F)\,\tilde{}\text{mod}\iso Z^c$mod. 
The equivalences $e_s$ for different $s$ are
(non-canonically) isomorphic.

{\it Proof.} For any $T(F)\,\tilde{}$-module $M$ the
action of $T\subset T(F)\,\tilde{}$ defines a
$\Gamma^\vee$-grading  $M = \oplus M^{\check{\gamma}}$ such
that
$T(F)^\gamma \,\tilde{}$ acts by operators of degree
$c(\gamma )$. The corresponding $\Gamma^\vee /c(\Gamma
)$-grading is compatible with the
$T(F)\,\tilde{}$-action; it is the decomposition of
$M$ by $Z^c$-isotypical components (one has
$\Gamma^\vee /c(\Gamma )=\Hom (Z^c ,\Bbb G_m )$).     

Set $M^s := \mathop\oplus\limits_{\chi \in
\Gamma^\vee /c(\Gamma )} M^{s(\chi )}$ and $e_s (M) :=
 (M^{s})^{T (\fm )} $ where
$T(\fm ):=
\Ker (T(O) \twoheadrightarrow T)$ is the unipotent radical
of $T(O)$. Then $T =T(O)/T(\fm )$ acts on a
$\chi$-isotypical component of $e_s (M)$ by the
characters $s(\chi )$, and $M$ is equal to
Ind$_{T(O)}^{T(F)\,\tilde{}} e_s (M)$.\footnote{It suffices
to show that $M^s =$ Ind$_{T(O)}^{T(F)^0 \,\tilde{}} e_s
(M)$. One has $T(F)^0 \,\tilde{} = T\times (T(F)^0
\,\tilde{}/T)$, and $T(F)^0 \,\tilde{}/T$-modules are the
same as the corresponding Lie algebra modules such that
Lie$T(\fm )$ acts by locally nilpotent operators. So our
assertion follows, say, from the usual Kashiwara's lemma. }
This induction is the functor inverse to
$e_s$.

Different $s$' differ by a map $\Gamma^\vee /c(\Gamma )\to
\Gamma \buildrel{c}\over\subset \Gamma^\vee$. An
isomorphism between the $e_s$'s is given by
a lifting of this map to $T(F)\,\tilde{}$ (here $\Gamma =
T(F)\,\tilde{}/T(F)^0 \,\tilde{}$).
\hfill$\square$

 \bigskip

\centerline{\bf \S2 Moving points} 

\medskip

{\bf 2.1} Let $X$ be a smooth curve. Then for every $x\in
X$ we have a local field $F_x$ of Laurent series at $x$.
The group ind-schemes $T(F_x )$ are fibers of a formally
smooth commutative group ind-scheme $T(F_X )$ over $X$.
For a $k$-algebra $R$ one has $T(F_X )(R)$ is the group of
pairs $(x, f^\gamma )$ where $x$ is an $R$-point of $X$
and $f^\gamma$ is a $T$-point with values in the
ring of functions on the punctured formal tubular
neighbourhood of the graph of
$x$. We have a group subscheme $T(O_X )\subset T(F_X )$
affine over $X$ and the quotient ind-scheme $Gr_X :=T(F_X
)/T(O_X )$ which is ind-finite, $Gr_X =\sqcup Gr_X^\gamma$.

Infinitesimal automorphisms of $X$ act on
$T(F_X )$ by transport of structure. Thus $T(F_X
)$ carries a canonical action of the Lie algebroid of
vector fields on $X$; for 
$\theta\in
\Theta (X)$ its action is denoted by Lie$_\theta$. Our
$T(F_X )$ also carries a canonical connection
$\nabla$ along $X$.  For $\theta \in\Theta (X)$ the
vertical vector field Lie$_\theta - \nabla_\theta$ on
$T(F_X )$ acts on every fiber $T(F_x )$ as the restriction
of
$\theta$ to $F_x$, see 1.5(iv). These three actions
preserve $T(O_X )$, so they also act on $Gr_X$.

 It is clear what a symmetric
Heisenberg extension of level $c$ in the present setting
is. Such an extension yields a superline $\lambda_{Gr_X} :=
T(F_X )\,\tilde{}/i(T(O_X )$ on $Gr_X$. Notice that the
restriction of $\lambda_{Gr_X}$ to $Gr_{X red}^{-\gamma}
=X$ is the superline $\lambda^\gamma_X$ from Lemma 1.7. The
latter Lemma remains true obviously modified, so symmetric
Heisenberg extensions form a
$\Gamma^\vee \otimes\mu_2$-gerbe on $X$ which is always
trivial (but has no canonical trivialization). As in
1.5(iv) we see that the above three actions of $\Theta
(X)$ on $T(F_X )$ lift canonically to any symmetric
Heisenberg extension $T(F_X )\,\tilde{}$.

{\bf 2.2} As every ind-schemes of jets and meromorphic
jets, $T(O_X )$ and
$T(F_X )$ have canonical factorization structures (see
[BD] 3.4). So for every $n\ge 1$ there is a canonical
ind-affine group ind-scheme
$T(F_{X^n})$ on $X^n$ equipped with a flat connection and
its affine group subscheme
$T(O_{X^n} )$ whose fibers
$T(F_{x_1 ,..,x_n )})$ over a
$k$-point $(x_1 ,..,x_n )\in X^n$ are products of fibers
$T(F_{x})$, $T(O_x )$ for $x\in \{ x_1 ,..,
x_n \}\subset X$. In particular, on $X\times X$ we have
$T(O_{X\times X})\subset T(F_{X\times X})$ whose pull-back
to the diagonal $\Delta : X\hra X\times X$ equals $T(O_X
)\subset T(F_X )$ and to the complement of the diagonal
$j:U\hra X\times X$ equals $T(O_X
)\times T(O_X )\subset T(F_X )\times T(F_X )$.

Recall the construction. First one has canonical
factorization commutative ring scheme and ind-scheme
$O_{X^\cdot}\subset F_{X^\cdot}$ equipped with connections.
Their
$R$-points can be described as follows. For $(x_1 ,.., x_n
)\in X^n (R)$ denote by $O_{(x_1 ,.., x_n )}$ the ring of
functions of the formal completion of $X\times \Spec R$ at
the union of graphs $Y_{x_i}$ of $x_i$'s, and by $F_{(x_1
,..,x_n )}$ the localization of $O_{(x_1 ,.., x_n )}$ with
respect to the equations of the divisors $Y_{x_i}$. These
$O_{(x_1 ,..,x_n )}\subset F_{(x_1 ,..,x_n )}$ are fibers
of
$O_{X^n}(R)\subset F_{X^n}(R)$ over $(x_1 ,..,x_n )$. The
connections and the factorization identifications are
obvious.

Now $T(O_{X^n})(R):= T(O_{X^n}(R))$ and $T(F_{X^n})(R):=
T(F_{X^n}(R))$.

The schemes $T(O_{X^n})$ are flat over $X^n$, and
ind-schemes $T(F_{X^n})$ ind-flat. 

{\bf 2.3} The Contou-Carr\`ere symbol satisfies the
factorization property, i.e., for every $n\ge 1$ we have a
pairing $\{
\, ,\, \}_n : F^\times_{X^n}\times F^\times_{X^n} \to
\CO_{X^n}^\times$ that satisfy the usual compatibilities.
Indeed, the definition of $\{\, ,\, \}$ from [BBE] \S2
immediately extends to the above situation; the
factorization compatibilities are clear.

A {\it factorization structure} on a symmetric Heisenberg
extension $T(F)\tilde{}_{X}$ is formed by symmetric
superextensions $T(F)\tilde{}_{X^n}$ of $T(F)_{X^n}$
connected by the usual factorization identifications (so
the restriction of  $T(F)\tilde{}_{X^n}$ to the fiber $\Pi
T(F_{x})$ as above is the Baer product of the symmetric
Heisenberg extensions
$T(F_x )\,\tilde{}\,$). The corresponding commutator
pairing is\footnote{Since
$T(F)\tilde{}_{X^n}$ is formally smooth and connected it
suffice to check this over the complement to the diagonal
divisor; now use the factorization property.} 
$\{ f_1^{\gamma_1}, f_2^{\gamma_2}\}^c_n = \{
f_1 ,f_2 \}_n^{-c(\gamma_1 ,\gamma_2 )}$. As follows from
Remark in 1.4, every symmetric Heisenberg extension admits
a unique factorization structure. The canonical flat
connection and the actions of $\Theta (X)$ on
$T(F_{X^n})$ lift canonically to $T(F_{X^n})\,\tilde{}$ in
the way compatible with the factorization structure.

The Lie algebra of $T(F_{X^n})\,\tilde{}$ is the symmetric
Heisenberg extension $\ft (F_{X^n} )\, \tilde{}$ of $\ft
(F_{X^n} )=\Gamma \otimes F_{X^n}$.

{\bf 2.4} The above factorization structures yield a
factorization structure on $Gr_X$
and the superline $\lambda$ on it. The 
action of the factorization group ind-scheme
$T(F_{X^n})$ on $Gr_{X^n}$ lifts to a (left) action of 
$T(F_{X^n})\,\tilde{}$ on $\lambda_{Gr_{X^n}}$.

Every $Gr_{X^n}$ is inductive limit of
subschemes finite and flat over $X^n$.

The factorization ind-scheme $Gr_X$ can be also understood
as the affine Grassmannian for the group $T$. Namely,
$R$-points of
$Gr_{X^n}$ are the same as triples $((x_1 ,.., x_n ), \CL
,\gamma )$ where
$(x_1 ,.. x_n )\in X^n (R)$, $\CL$ is a $T$-torsor over
$X\times \Spec R$, $\gamma$ is a trivialization of $\CL$
over the complement to the graphs of $x_i$. 

{\it Remarks.} (i) The distinguished horizontal section 1
of
$Gr_X$ is a unit for the factorization structure, i.e.,
for any (local) section $\phi$ of $Gr_X$ the section
$\phi \times 1$ of $Gr_{X\times X}$ defined a priori on
$U$ extends to a section over $X\times X$ whose
restriction to the diagonal equals $\phi$. The superline
$\lambda_{Gr_X}$ is canonically trivialized over 1, and
this trivialization is compatible with the factorization
structure on $\lambda_{Gr_X}$.

(ii) The
section 1 of $T(O_X )$ is {\it not} a unit for the
factorization structure on $T(O_X )$.

{\bf 2.5} (the reader can skip it). The affine Grassmannian
description makes it clear that for every $n\ge 2$ there
is a natural morphism $\pi^n : Gr_{X}^n \to Gr_{X^n}$ which
equals the factorization isomorphism over the complement
to the diagonal divisor in $X^n$. Namely, $\pi^n$ sends
$((x_1 ,\CL_1 ,\gamma_1 ),..,(x_n ,\CL_n ,\gamma_n ))$ to
$((x_1 ,..,x_n ),\CL_1 \cdot ..\cdot\CL_n ,\gamma_1 \cdot
..\cdot\gamma_n )$ where $\cdot$ is
 the Baer product. 

{ Lemma.} The factorization isomorphism yields an
isomorphism of superlines 
$$\mu_n :\lambda_{Gr_X^{\gamma_1}}\boxtimes
..\boxtimes\lambda_{Gr_{X}^{\gamma_n}}\iso (\pi^{n*}
\lambda_{Gr_{X^n}})(-
\mathop\sum\limits_{i<j} c(\gamma_i ,\gamma_j
)\Delta_{ij} )$$
on $Gr_{X}^{\gamma_1}\times
..\times Gr_X^{\gamma_n}$ (here $\Delta_{ij}$ are the
diagonal divisors $x_i =x_j$). The restriction of
$\mu_2$ to
$Gr_{X red}^{\gamma_1} \mathop\times\limits_X  Gr_{X
red}^{\gamma_2} =X$ is the isomorphism $\mu_{\gamma_1
,\gamma_2}$ from Lemma 1.7.
 
{\it Proof.} It suffices to consider the case of $n=2$
and restrict our picture to the reduced schemes. One has
$Gr_{X red}^{\gamma_i} \iso X$, so we have superlines
$\lambda^i_X :=\lambda_{Gr_{X red}^{\gamma_i}}$ on $X$ and
$\lambda_{X\times X}$ on $X\times X$ defined as the
restriction of
$\pi^{2*} \lambda_{Gr_{X\times X}}$ to 
$Gr_{X red}^{\gamma_1}\times Gr_{X red}^{\gamma_2}$.
The factorization  yields an isomorphism
$\phi :
\lambda^1_X
\boxtimes\lambda^2_X \iso \lambda_{X\times X}(?\Delta )$; 
we want to show that $? = -c(\gamma_1 ,\gamma_2 )$ and
compute the restriction of $\phi$ to the diagonal.

Let $t$ be a local coordinate on $X$, $x$ the
corresponding coordinate on  another copy of $X$, $x_1
,x_2$ the  local coordinates on $X\times X$. Recall that
$\lambda_{Gr_{X }^{\gamma}}$ is formed by sections of
$T(F_X )\,\tilde{}$ over $T(F_X)_{red}^\gamma$ which are
right invariant with respect to the $i(T(O_X ))$-action.
Therefore 
$\lambda^i_X$ identifies with the pull-back of $T(F_X
)\,\tilde{}$ by the section $x\mapsto (t-x)^{\gamma_i}$
of $T(F_X )^{\gamma_i}$. Similarly, the superlines $p_i^*
\lambda^i_X$ on $X\times X$  identify with the pull-back
of $T(F_{X\times X} )\,\tilde{}$ by the sections $(x_1
,x_2 )\mapsto (t-x_i )^{\gamma_i}$ of $T(F_{X\times X})$,
and
$\lambda_{X\times X}$ is the pull-back by the product of
these sections. These identifications together with the
product in $T(F_{X\times X} )\,\tilde{}$ provide an
isomorphism
$\psi =\psi_t :
\lambda^1_X \boxtimes\lambda^2_X \iso \lambda_{X\times X}$.

It remains to compute the function $\phi /\psi$. Let
$\ell^i$ be invertible (local) sections of $\lambda^i_X$ 
considered as sections of the
pull-back of $T(F_{X })\,\tilde{}$ by $x\mapsto
(t-x )^{\gamma_i}$. Our $\ell^i$ yields an invertible
section $\hat{\ell}^i$ of  
$ p_i^* \lambda^i_X$ identified with the
pull-back of $T(F_{X \times X})\,\tilde{}$ by $x\mapsto
(t-x_i )^{\gamma_i}$. 
The product $\hat{\ell}^1 \cdot\hat{\ell}^2$ is a
section of the pull-back of $T(F_{X \times X})\,\tilde{}$
by $x\mapsto (t-x_1)^{\gamma_1}(t-x_2 )^{\gamma_2}$. From the point of view of the
factorization isomorphism outside of the diagonal our
sections are
$\hat{\ell}^1_{(x_1 ,x_2 )} = \ell^1_{x_1} i(t-x_1
)^{\gamma_1}_{x_2}$, $\hat{\ell}^2_{(x_1 ,x_2 )}= i(t-x_2
)^{\gamma_2}_{x_1}\ell^2_{x_2}$, and $(\hat{\ell}^1
\cdot\hat{\ell}^2 )_{(x_1 ,x_2 )}= (\ell^1 i(t-x_2
)^{\gamma_2})_{x_1} (i(t-x_1 )^{\gamma_1}\ell^2 )_{x_2}$.
To compute
$\phi /\psi$ we have to extend
$\ell^1 \boxtimes \ell^2$ and $\hat{\ell}^1
\cdot\hat{\ell}^2$ to sections of $T(F_{X\times X
})\,\tilde{}$ invariant with respect to the right
$iT(O_{X\times X})$-translations and take the ratio.
Thus $\phi /\psi = \ell^2_{x_2}/\text{Ad}_{i(t-x_1
)^{\gamma_1}}\ell^2_{x_2}= (x_2 -x_1 )^{c(\gamma_1
,\gamma_2 )}$. This implies Lemma.   \hfill$\square$

 \bigskip

\centerline{\bf \S3 The vertex counterpart} 

\medskip

{\bf 3.1} We  use the framework of [BD], so for us a
vertex algebra
$A$ is the same as a universal chiral (super)algebra. We
consider it as a chiral algebra $A_O$ over Spec$(O)$
equipped with an action of Aut$^{1/2}(O)$ (so $A=A_O
/\fm_x A_O$). Sometimes we prefer to consider the
corresponding chiral algebra $A_X$ on a curve
$X$. Similarly, ``vertex Lie algebras" from [FBZ]
or ``conformal algebras" from [K] are the same as universal
Lie$^*$ algebras. For such a fellow
$L$ we have the corresponding Lie algebras $h_F (L)\supset
h_O (L)$ of 
$\partial_t$-coinvariants on $L_F \supset L_O$. We always
have a distinguished central element $1_L \in L$; the
 enveloping chiral algebra is denoted by
$U(L)$.\footnote{Our $U(L)$ is the plain chiral envelope
modulo the relation $1_L =1$.}

We denote by $\ft\tilde{_\CD}$ the symmetric Heisenberg
Lie$^*$ algebra of level $c$. Then
$h_F (\ft\tilde{_\CD})$ is the symmetric Heisenberg Lie
algebra $\ft (F)\,\tilde{}$ from 1.5(ii).

{\bf 3.2 Definition.} A {\it symmetric lattice vertex
algebra of level $c$} is a $\Gamma$-graded vertex algebra
$A=\oplus A^\gamma$  equipped with a morphism of
Lie$^*$ algebras
$\alpha : \ft\tilde{_\CD}
\to A^0$ and an automorphism $\sigma$
such that:

(a) $\alpha$ is compatible with $\sigma$, sends $1\in
\ft\tilde{_\CD}$ to $1\in A$, and
yields $U(\ft\tilde{_\CD})\iso A^0$. 

(b) Each $A^\gamma$ has pairity
$c(\gamma ,\gamma )$ mod2. As a $\ft
(F)\,\tilde{}$-module $A^\gamma$ is isomorphic to the
representation Ind$c(\gamma )$ induced by the character
$c(\gamma )\in \Gamma^\vee =\Hom (T(O),\Bbb G_m )$ of
the Lie subalgebra $\ft (O)\buildrel{i}\over\hra
\ft (F)\,\tilde{}$.

(c)  The chiral product on every pair of
components $A^{\gamma_1}$, $A^{\gamma_2}$ does not vanish.

(d) One has $\sigma^2 =$ id$_A$ and $\sigma (A^\gamma
)=A^{-\gamma}$.

{\it Example.} The algebra $A$ of functions on the scheme
of $T^\vee$-valued jets is a symmetric Heisenberg vertex
algebra of level 0. Here the $\Gamma$-grading comes from
the translation action of $T^\vee$, $\sigma$ comes from
the involution $\tau \mapsto \tau^{-1}$ of $T^\vee$, and
$\alpha$ from the embedding $\ft \otimes\omega_x^{-1}
=\Gamma \otimes\theta_x^{-1}\hra A$, $\gamma
\otimes\partial_x \mapsto e^{-\gamma}\partial_x e^\gamma
$ where $e^\gamma$ denotes the character of $T^\vee$
corresponding to $\gamma \in\Gamma$. 

{\it Remark.} If $c$ is non-degenerate then representations
Ind$c(\gamma )$ are irreducible and pairwise
non-isomorphic, so the datum of $\Gamma$-grading is
redundant.

 \medskip

{\bf 3.3} Let $T (F)\, \tilde{}$ be a symmetric
Heisenberg extension of level $c$. Denote by $A$ the
representation of $T(F)\,\tilde{}$ induced from the
trivial representation of
$T(O)\buildrel{i}\over\hra T(F)\,\tilde{}$.

{\bf Proposition.} $A$ is a
symmetric lattice vertex algebra in a natural way. 

{\it Proof.} A comment on the definition of $A$. As in 2.1
set $Gr := T(F)/T(O)$; this is an ind-finite
ind-scheme whose connected components are labeled by
$\gamma \in \Gamma$. Denote by
$\lambda_{Gr}$ the superline $T(F)\,\tilde{} /i(T(O))$
over $Gr$; it is a $T(F)\, \tilde{}$-equivariant
superline in the obvious way. Let $\delta_{Gr}$ be the
dualizing $\CO$-module on $Gr$. Equivalently,
$\delta_{Gr}$ is the cofree $T(F)$-equivariant
$\CO$-module on
$Gr$ such that
$i_1^! \delta_{Gr} =k$. Set
$A^\gamma :=
\Gamma (Gr^{-\gamma }, \lambda_{Gr}\otimes\delta_{Gr})$ and
$A:=\oplus A^\gamma $.

We see  that $A$ is
$\Gamma$-graded  and is equipped with an action of
$\sigma$ and Aut$^{1/2}(O)$ (they act on
$T(F)\,\tilde{}\,$ preserving $T(O)$). In particular, the
Lie algebra of vector fields $\Theta (O)$ acts on $A$.
Each $A^\gamma$ contains the superline $\lambda^\gamma :=
\Gamma (Gr^{-\gamma}_{red},\lambda_{Gr})$ which equals
$\lambda^\gamma$ from 1.7.

For a curve $X$ and a symmetric Heisenberg extension
$T(F_X )\,\tilde{}$ on $X$ we get a $\CD_X$-module $A_X$.
A vertex algebra structure on $A$ amounts to a natural
chiral algebra structure on $A_X$. Now the
factorization structure on $Gr$ and
$\lambda_{Gr}$ as defined in 2.4 define a factorization
structure on $A_X$ which amounts to a chiral algebra
structure (see [BD] 3.4). Thus $A$ is a vertex algebra
with unit $1_A =\lambda^0$.

The group ind-schemes $T(F_{X^n})\,\tilde{}$
act on $\CO_{X^n}$-modules $A_{X^n}$ (see 2.4) in a way
compatible with the $\CD$-module and factorization
structures.

Notice that $A^0$ equals the chiral
envelope of $\ft\tilde{_\CD}$,\footnote{ See [BD] 
3.7.5--3.7.11, 3.7.20.} so 3.2(a) is valid. Properties
3.2(c),(d) are clear. 
 To see 3.2(b) note that $\ft
(F)\,\tilde{}$ acts on $A$ in two ways: first as the Lie
algebra of
$T(F)\,\tilde{}$ and second as
$h_F (\ft\tilde{_\CD}) \subset h_F (A)$ by the
adjoint chiral algebra action. The two actions
coincide. For the first action it is clear that 
$\ft (O) \subset \ft (F)\,\tilde{}$ acts on the superline 
$\lambda^\gamma \subset A^\gamma$ by the
character $c(\gamma )$, and $A^\gamma$ is the
corresponding induced $\ft (F)\,\tilde{}$-module. 
\hfill$\square$

\medskip

{\bf 3.4 }  Here is another  description of the chiral
product on $A_X$. The action of $T(F_X
)\,\tilde{}$ on $A_X$ yields, via the factorization
structure on $T(F_X )\,\tilde{}$, an action of
$T(F_{X^2})\,\tilde{}$  on $j_* j^* (A_X
\boxtimes A_X )$ and $\Delta_* A_X$ compatible with the
connections.\footnote{Notice that for a $\CD_X$-module
$M$ a
$T(F_X )\,\tilde{}$-action on $\Delta_* M$  compatible with
the connections amounts to a $\Delta^* T(F_X
)\,\tilde{}$-action on $M$ compatible with the
connections.} Now the embedding $1_A : \CO_X \hra A_X$
identifies $A_X$ with the $T(F_X )\,\tilde{}$-module
induced from the trivial $T(O_X )$-module $\CO_X$.
So $1_A^{\boxtimes 2} : j_* j^*
\CO_{X\times X} \hra j_* j^* (A_X
\boxtimes A_X )$ identifies $j_* j^* (A_X
\boxtimes A_X )$ with the $T(F_{X^2})\,\tilde{}$-module
induced from the trivial $T(O_{X^2})$-module $j_* j^*
\CO_{X\times X}$. Therefore the canonical morphism $j_*
j^* \CO_{X\times X}\to \Delta_* \CO_X$ extends canonically
to a morphism of $T(F_{X^2})\,\tilde{}$-modules 
$$\mu_A : j_* j^* (A_X
\boxtimes A_X )\to\Delta_* A_X . \tag 3.4.1$$

{\bf Lemma.} The chiral product on $A_X$ equals 
$\mu_A$.

{\it Proof.} Recall (see 3.3) that 
$T(F_{X^2})\,\tilde{}$ acts naturally on
$A_{X^2}$, hence on $j_* j^* A_{X^2}$ and $j_* j^*
A_{X^2}/A_{X^2}$. The factorization isomorphisms $j_* j^*
A_{X^2}\iso j_* j^* (A_X \boxtimes A_X )$, $j_* j^*
A_{X^2}/A_{X^2}\iso \Delta_* A_X$ are compatible with the 
$T(F_{X^2})\,\tilde{}$-actions. Thus the chiral product on
$A_X$ is compatible with the
$T(F_{X^2})\,\tilde{}$-action. On $1_A$ it equals the canonical morphism $j_*
j^* \CO_{X\times X}\to \Delta_* \CO_X$. We are done
since these two properties characterize $\mu_A$.
\hfill$\square$

\medskip

{\bf 3.5} Let us deduce the usual explicit formula for 
vertex operators. For a point $x\in X$ consider the
operator product 
$\circ : A_{O_x}
\otimes A_{x}\to F_x
\hat\otimes A_x$, $\phi\otimes a \mapsto \phi\circ a$,
with fixed second variable equal to $x$. An
invertible section
$\ell\in
\lambda^\gamma_{O_x}$ defines then the
vertex operator
$V^\gamma_\ell : A_x \to F_x
\hat\otimes A_x$, $V^\gamma_\ell (a):=  \ell\circ a$,
 which we want to compute.  
 
Denote by $F_{O_x}$ the localization of
$O_x\hat\otimes O_x$ with respect to an equation of the
diagonal. One has obvious morphisms $ F_{O_x} \to F_x
\hat\otimes O_x := \limleft F_x \otimes (O_x /\fm^n_x )$
and
$ F_{O_x} \to O_x \hat\otimes F_x := \limleft (O_x
/\fm_x^n )\otimes F_x$ which we denote by $f \mapsto f_+
, f_-$. 

Let
$g\in T(F_{O_x} )$ be any element such that its order of
pole along the diagonal divisor equals $\gamma$; such $g$
is defined up to multiplication by an element of $T(O_x
\hat\otimes O_x )$. So
 $g$ is a section of $T(F_X )^{-\gamma}$ over Spec$\,
O_x$. Thus $\lambda_{O_x}^\gamma$ equals the
pull-back of $T(F_X )\,\tilde{}$ by $g$, so our $\ell$ 
defines a lifting
$\tilde{g}:=\ell (g) \in T(F_X
)\,\tilde{}$  of $g$.

Two points of view on $g$:

(i) We have  $g_+ \in T(F_x )(O\hat{_x} ):= \limleft
T(F_x )(O_x /\fm_x^n )$. In other words,
$g_+$ is the restriction of $g$ to the formal
neighbourhood of $x$ identified with an $O\hat{_x}$-point
of
$T(F_x )$ by means of the connection on $T(F_X )$.
Similarly, we have $\tilde{g}_+ \in T(F_x )\tilde{}\,
(O\hat{_x})$.

(ii) We have $g_- \in T(O_x )(F_x )$ hence
 $\tilde{g}_- := \tilde{i}(g_- )\in T(O_x )\tilde{}
(F_x )$.

Consider $\tilde{g}$ as a section of the
restriction of $T(F_{X^2})\,\tilde{}$ to $
\Spec O_x \times \{ x\}
\subset X\times X$. Then $\tilde{g}$ acts on $
A_{F_x}\otimes A_x$ (which is the restriction of $j_* j^*
(A_X \boxtimes A_X )$) as $\tilde{g} \otimes
\tilde{g}_-$, and on $ (F_x /O_x )\otimes A_x$
(which is the restriction of $\Delta_* A_X$) as
$\tilde{g}_+$. Since the chiral product $\mu_A$, so the
ope product, commutes with the action of $\tilde{g}\in
T(F_{X^2})\,\tilde{}$ (see 3.4) we get, applying them to
$1_A \otimes \tilde{g}_{-}^{-1} a
 \in  F_x
\cdot1_A \otimes A_x
\subset  A_{F_x}\otimes A_x$, the formula 
$$V^\gamma_{\ell} = \tilde{g}_+
\tilde{g}_-^{-1}. \tag 3.5.1$$

This is, indeed, the usual formula. To see this choose a
coordinate $t$ at $x$; let $z$ be the same coordinate on
another copy of $X$. Pick any non-zero $\ell_x \in
\lambda_x^\gamma$; let $\ell \in \lambda^\gamma_{O_x}$
be the corresponding translation invariant\footnote{i.e.,
Lie$_{\partial_z}$-invariant.} section. Take
$g:= (t-z)^{-\gamma}\in T(F_{O_x} )$; as above, $\ell$
yields its translation invariant lifting $\tilde{g}:=\ell
(g)$ which at $z=0$ takes value
$\tilde{t}^{-\gamma}:=\ell_x (t^{-\gamma})\in T(F_x
)\,\tilde{}$. Set
$h_a^\gamma :=\tilde{i}(t^a
\cdot \gamma ) \in \ft (F_x )\,\tilde{}$ (see 1.5, 1.6).  
Then
$\tilde{g}_+ (z) = \exp
(\mathop\sum\limits_{n>0}h_{-n}^\gamma
\frac{z^n}{n})  \tilde{t}^{-\gamma}$, $ \tilde{g}_-
(z)= \exp 
(\mathop\sum\limits_{n>0}h_{n}^\gamma
\frac{z^{-n}}{n})(-z)^{-h_0^\gamma}$, and
$$V^\gamma_\ell = \exp
(\mathop\sum\limits_{n>0}h_{-n}^\gamma
\frac{z^n}{n})  \tilde{t}^{-\gamma}\exp 
(-\mathop\sum\limits_{n>0}h_{n}^\gamma
\frac{z^{-n}}{n})(-z)^{h_0^\gamma}. \tag 3.5.2$$

{\bf 3.6 Remarks.}
(a) For a coordinate
$t$ as above set
$L_a := t^{a+1}\partial_t \in \Theta (F)$.
Now, by 1.6, for $\ell \in \lambda^\gamma \subset A^\gamma$
one has
$$L_{-1}\ell = -h_{-1}^\gamma \ell, \quad L_0 \ell =
-\frac{1}{2}c (\gamma ,\gamma )\ell ,  \quad L_{\ge 1}\ell
=0. \tag 3.6.1$$
These formulas, together with the compatibility of the
$\Theta (F)$-action with the $\ft (F)\, \tilde{}$-action
recover the $\Theta (F)$-action on $A$ via the 
$\ft (F)\, \tilde{}$-action. 

The first of the above
equations implies that the operator $V^\gamma_\ell$ from
(3.5.2) satifies the following differential
equation:\footnote{Which is clear also from the formula
(3.5.2) itself.} $$z\partial_z V^\gamma_\ell =
\mathop\sum\limits_{n>0}h_{-n}^\gamma  V^\gamma_\ell z^n +
 \mathop\sum\limits_{n\ge 0} V^\gamma_\ell h_{n}^\gamma
z^{-n} . \tag 3.6.2$$

(b) Let $x_1$, $x_2$ be the coordinates on $X\times X$ that
correspond to $t$. As follows from (3.5.2), for every
$\gamma_1
,\gamma_2
\in
\Gamma$ 
and
$\ell_i \in \lambda^{\gamma_i }$ the operator product
$\ell_1
\circ\ell_2$ belongs to $(x_1 -x_2 )^{c(\gamma_1 ,\gamma_2
)} A^{\gamma_1 +\gamma_2}_{x_2}$. The top coefficient
lies, in fact, in
$ (x_1 -x_2)^{c(\gamma_1 ,\gamma_2 )}\lambda^{\gamma_1
+\gamma_2}_{x_2}$ and it equals $\mu (\ell_1 ,\ell_2 )$
where $\mu$ was defined in 1.7.\footnote{Precisely, the
top coefficient is
$ (x_1 -x_2)^{c(\gamma_1 ,\gamma_2
)}\mu (\ell_1 ,\ell_2 )dt^{-c(\gamma_1 ,\gamma_2
) }$.}

(c)\footnote{Need not for the rest of the text.} The chiral
algebra structure on
$A$ can be described
 in terms of the
$\ft (F)\, \tilde{}$-action as follows. It suffices
to recover   $A_{X^2}^{\gamma_1
\gamma_2}:= A_{X^2}\cap j_* j^* A_X^{\gamma_1}
\boxtimes A_X^{\gamma_2}\subset j_* j^* A_X^{\gamma_1}
\boxtimes A_X^{\gamma_2}$ together with  identifications
$\Delta^* A_{X^2}^{\gamma_1 \gamma_2}\iso A_X^{\gamma_1
+\gamma_2}$.
 Consider the sheaves of Lie $\CO_{X^2}$-algebras
$\ft (O_{X^2})\subset \ft (F_{X^2}) \subset j_* j^* \ft
(F_X )\times \ft (F_X )$ equipped with a connection
$\nabla$ and the symmetric Heisenberg extension $\ft
(F_{X^2})\,\tilde{}$ (see 2.2).  

{\bf Lemma.} 
 $A_{X^2}^{\gamma_1
\gamma_2}$ is the $\ft (F_{X^2})\,\tilde{}$-submodule of 
$j_* j^* A_X^{\gamma_1}
\boxtimes A_X^{\gamma_2}$ generated by the superline 
$\lambda_{X^2}^{\gamma_1 \gamma_2}:= \lambda^{\gamma_1}_{
X}\boxtimes
\lambda^{\gamma_2}_{ X} (c(\gamma_1 ,\gamma_2 )\Delta )$.
The isomorphism
$\Delta^* A_{X^2}^{\gamma_1 \gamma_2}\iso A_X^{\gamma_1
+\gamma_2}$ is the morphism of $\ft
(F_{X^2})\,\tilde{}$-modules induced by the morphism
$\mu_{\gamma_1 \gamma_2}
:\Delta^*
\lambda_{X^2}^{\gamma_1
\gamma_2}
=\lambda^{\gamma_1}_X \otimes\lambda^{\gamma_2}_X \otimes
\omega_X^{\otimes -c(\gamma_1 ,\gamma_2 )}\iso
\lambda_X^{\gamma_1 +\gamma_2}$ between the generators
(see 1.7).

{\it Proof.} By 2.5 $A_{X^2}^{\gamma_1
\gamma_2}$ contains the superline  $\lambda^{\gamma_1
\gamma_2}_{ X^2}$ hence the $\ft
(F_{X^2})\,\tilde{}$-submodule 
$A_{X^2}^{'\gamma_1 \gamma_2}$ generated by it.
 $A_{X^2}^{'\gamma_1
\gamma_2}$ equals the
$\ft (F_{X^2})\,\tilde{}$-module induced from the $\ft
(O_{X^2})$-module 
$\lambda^{\gamma_1\gamma_2}_{ X^2}$ on which $\ft
(O_{X^2})$ acts by the character $\phi (t)_{x_1
,x_2}\mapsto
\gamma_1
\phi (x_1 )_{x_1 ,x_2}+
\gamma_2 \phi (x_2 )_{x_1 ,x_2}$. So  $\Delta^*
A_{X^2}^{'\gamma_1 \gamma_2}$ is the
$\Delta^* \ft (F_{X^2})\,\tilde{} =\ft (F_X
)\,\tilde{}$-module induced from the $\ft (O_X )$-module
$\Delta^* \lambda^{\gamma_1
\gamma_2}_{ X^2}
=\lambda_X^{\gamma_1}\otimes\lambda_X^{\gamma_2}\otimes
\omega_X^{\otimes -c(\gamma_1 ,\gamma_2 )}$ on which $\ft
(O_X )$ acts by 
$\gamma_1 +\gamma_2$. 

To show that $A_{X^2}^{'\gamma_1
\gamma_2}=A_{X^2}^{\gamma_1
\gamma_2}$ it suffices to check that
  $A_{X^2}^{'\gamma_1
\gamma_2}$ is a $\CD_{X^2}$-submodule (since both $A_{X^2}^{'\gamma_1
\gamma_2}$, $A_{X^2}^{\gamma_1
\gamma_2}$ have the same pull-back to
the diagonal) which in turn amounts to the fact
 that the action of vector fields sends
 $\lambda_{X^2}^{\gamma_1 \gamma_2}$ to
$A^{'\gamma_1 \gamma_2}_{X^2}$. Let $t$ be a local
coordinate on $X$,
$x_1$,
$x_2$ the corresponding coordinates on $X^2$. For
translation invariant\footnote{Here ``translation" is the
action of Lie$\partial_t$.} sections
$\ell_i
\in \lambda_X^{\gamma_i}$
one has
$\partial_{x_1} ((x_1 -x_2 )^{-c(\gamma_1 ,\gamma_2
)}\ell_1
\boxtimes\ell_2 )= -c(\gamma_1 ,\gamma_2 )(x_1
-x_2 )^{-c(\gamma_1 ,\gamma_2 )-1}\ell_1
\boxtimes\ell_2 +(x_1 -x_2)^{-c(\gamma_1
,\gamma_2 )}(h_{-1}^\gamma \ell_1 )
\boxtimes\ell_2 .$ Consider
$\phi :=\tilde{i}((t-x_1 )^{-1}\cdot\gamma )\in \ft
(F_{X^2})\,\tilde{}$; one has $\phi (\ell_1 \boxtimes\ell_2
)= (h_{-1}^\gamma \ell_1 )\boxtimes \ell_2 - c(\gamma_1
,\gamma_2 )(x_2 -x_1 )^{-1}\ell_1
\boxtimes \ell_2$. Therefore  $\partial_{x_1}
((x_1 -x_2 )^{c(\gamma_1 ,\gamma_2 )}\ell_1
\boxtimes\ell_2 )= \phi (((x_1 -x_2 )^{c(\gamma_1 ,\gamma_2 )}\ell_1
\boxtimes\ell_2 )$.  

The last statement of Lemma follows from 2.5.
\hfill$\square$

\medskip

{\bf 3.7} 
Symmetric lattice vertex algebras
of level
$c$ form a groupoid
$\CV s^c$. In 3.3 we have defined a functor $\CH
s^c \to \CV s^c$, $T(F)\,\tilde{} \mapsto A = A(
T(F)\,\tilde{} )$. 

{\bf Lemma.} If $c$ is non-degenerate then our functor
$\CH
s^c \to \CV s^c$ is an equivalence.

{\it Proof.}  Let $A$ be a symmetric lattice vertex algebra
of level $c$. Set $\lambda^\gamma := (A^\gamma )^{T(\fm
)}$; this is a superline on which $T(O)$ acts by the
character $c(\gamma )$, and $A^\gamma$ is the
corresponding induced $\ft (F)\,\tilde{}$-module. 

(i) Let us
check that formulas (3.6.1) remain valid. Set $\ft_a :=
\tilde{i}(\ft t^a )\subset \ft (F)\,\tilde{}$; for $a\neq
0$ this is the $a$-eigenspace for the $L_0$-action on $\ft
(F)\,\tilde{}$. Take a non-zero
$\ell
\in
\lambda^\gamma$. It is clear that
$L_{\ge 1}\ell =0$ and
$L_0 \ell =w\ell$ for some $w\in k$. Since
$[L_{-1},L_0 ]=L_{-1}$ we see that $L_{-1}\ell$
is an eigenvector of
$L_0$ of eigenvalue $w-1$, hence $L_{-1}\ell = \phi \ell$
for some $\phi \in \ft_{-1}$. For every $\psi \in \ft_{1}$
one has $\psi L_{-1}\ell = -(L_{-1}(\psi))\ell = [\psi
,\phi ]\ell$. Since $c$ is non-degenerate this implies
that $L_{-1}\ell = -h_{-1}^\gamma \ell$. Now, since $L_0 =
\frac{1}{2}[L_{-1}, L_2 ]$ and $L_2 (h^\gamma_{-1})=-
h^\gamma_0$, one has
$L_0
\ell = \frac{1}{2}L_2
h_{-1}^\gamma \ell =-\frac{1}{2}h^\gamma_0 \ell =
-\frac{1}{2}c(\gamma ,\gamma )\ell$, and we are done.

(ii) Let us check that 3.6(b) remains valid.
 Since $\circ$ is compatible with the Heisenberg
action we see that that the top non-zero coefficient of
$\ell_1 \circ
\ell_2$ belongs to $\lambda^{\gamma_1 +\gamma_2}$.
By (i) we know how $L_0$ acts on $\lambda^\gamma$,
therefore the compatibility of $\circ$ with the
$L_0$-action implies that this coefficient has degree
$-\frac{1}{2}c(\gamma_1 ,\gamma_1) -\frac{1}{2}c(\gamma_2
,\gamma_2 ) +\frac{1}{2}c(\gamma_1 +\gamma_2 ,\gamma_1
+\gamma_2 ) =c(\gamma_1 ,\gamma_2) $.

(iii) The above top non-zero coefficients form a system of
isomorphisms $\mu :
\lambda^{\gamma_1}\otimes\lambda^{\gamma_2} \iso
\lambda^{\gamma_1 +\gamma_2}\otimes\omega_x^{c(\gamma_1
,\gamma_2 )}$. Together with identifications $\sigma
:\lambda^\gamma \iso \lambda^{-\gamma}$ they form the
datum of 1.7, so we get the 
symmetric Heisenberg extension $T(F)\,\tilde{}$.

(iv) It follows directly from the construction of 3.3
that the functor $A \mapsto T(F)\,\tilde{}$ is left
inverse to the functor from 3.3. Let us show that it is 
 right inverse as well. Let $A$ be as above,
$T(F)\,\tilde{}$ the Heisenberg extension just defined,
and let $A'$ be the lattice
vertex algebra defined by $T(F)\, \tilde{}$ as in 3.3. Then
there is a canonical isomorphism $A \iso A'$ of $\ft
(F)\,\tilde{}$-modules which is identity map on the
generators $\lambda^\gamma$. We have to show that this
isomorphism is compatible with the vertex algebra
structures (the compatibility with the $\Gamma$-gradings
and $\sigma$ is clear).

We already know that $A\iso A'$ is compatible with the
$\Theta (O)$-action. Denote by $\circ'$ the ope
product on $A$ that comes from the one of $A'$. To show
that $\circ =\circ'$ it suffices to prove that
$\ell_1 \circ \ell_2 =\ell_1 \circ' \ell_2 \in  A((z))$ for
$\ell_i \in \lambda^{\gamma_i}$. We know that these formal
power series have the same order of pole, the same leading
coefficient, and satisfy the same differential equation
$z\partial_z p(z) = \mathop\sum\limits_{n>0}
h^{\gamma_1}_{-n}p(z)z^n + c(\gamma_1 ,\gamma_2 )p(z)$
which arises from the first equality in (3.6.1). Hence
they are equal, q.e.d.
\hfill$\square$ 

\medskip

{\bf 3.8 Representations.} For $x\in X$ we have
the category of $T(F_x )\,\tilde{}$-modules $T(F_x
)\,\tilde{}$-mod (see 1.9) and the category $\CM_x (A )$
of $A_X$-modules supported at $x$.

{\bf Lemma.} There is a canonical fully faithful
embedding
$ T(F_x )\,\tilde{}\text{-mod}\to \CM_x (A). $

{\it Proof.}  The construction, parallel to 3.4, goes as
follows.  Notation:
$j_x :U_x := X\smallsetminus \{ x\} \hra X$, $i_x : \{ x\}
\hra X$. 

Let $O_{X,x} \subset F'_{X,x}\subset F_{X,x}$ be the
commutative ring (ind)-schemes over $X$ with fibers over
$y\in X(R)$ equal, respectively, to $O_{(y,x)}$, $F'_y$:=
the localization of $O_{(y,x)}$ with respect to the
equation of $x$, and $F_{(y,x)}$ (see 2.2). They carry an
obvious canonical connection. There is an obvious
projection $F'_{X,x} \to (F_x)_X$ and the factorization
identifications
$i_x^* F'_{X,x}=i_x^* F_{X,x}=F_x$,
$j_x^* F'_{X,x}= j_x^* O_X \times F_x$, $j_x^* F_{X,x}=
j_x^* F_X \times F_x$.

We have the corresponding group ind-schemes $T(O_{X,x})
\subset T(F'_{X,x})\subset T( F_{X,x})$. The latter one is
the restriction of $T(F_{X^2 })$
to $X\times \{ x\} \subset X\times X$, so we have the
corresponding superextension $T(F_{X,x})\,\tilde{}$. The
restriction of $T(F_{X,x})\,\tilde{}$ to $x\in X$ equals
$T(F_x )\,\tilde{}$.   The pull-back
$T(F'_{X,x})\,\tilde{}$ of $T(F_x )\,\tilde{}$ by the
projection $T(F'_{X,x})\to T(F_x )$ coincides with the
restriction of
$T(F_{X,x})\,\tilde{}$ to
$T(F'_{X,x})$.\footnote{To see this use e.g. Remark in
1.4.} 

 Let $M_x $ be a $T(F_x )\,\tilde{}$-module. We will
construct a canonical chiral action of $A$ on $i_{x*}M_x$
which defines (3.6.1). 
$T(F'_{X,x})\,\tilde{}$ acts on $j_{x*}j_x^* \CO_X \otimes
M_x$ via $T(F'_{X,x})\,\tilde{}\to T(F_x )\,\tilde{}$. The
corresponding induced $T(F_{X,x})\,\tilde{}$-module equals 
$j_{x*}j_x^* A_X \otimes M$. Notice that
$T(F_{X,x})\,\tilde{}$ acts on the $\CD_X$-module
$i_{x*}M_x$ in a way compatible with connections,\footnote{
Such action amounts to an $i_x^*
T(F_{X,x})\,\tilde{}$-action which we have since $i_x^*
T(F_{X,x})\,\tilde{} =T(F_x )\,\tilde{}$.} so the obvious
morphism $j_{x*}j_x^* \CO_X \otimes
M_x \twoheadrightarrow i_{x*} M_x$ extends to a morphism
of $T(F_{X,x})\,\tilde{}$-modules $\mu_M :j_{x*}j_x^* A_X
\otimes M_x \to i_{x*}M_x$. This is the promised chiral
$A_X$-action on $i_{x*}M_x$.

It is clear that the chiral action of the Lie Heisenberg
$\ft\tilde{_\CD}$ on $i_{x*}M_x$ coincides with the $\ft
(F_x )\,\tilde{}$-action. Repeating the argument of 3.5 we
see that
$\lambda^\gamma_X \subset A^\gamma_X$ acts on $i_{x*}M_x$
according to formula (3.5.1). These operators determine
completely the action of $T(F_x )\,\tilde{}$ on $M_x$, so
our functor is fully faithful.   \hfill$\square$

\medskip

{\bf 3.9 Proposition.} If $c$ is non-degenerate then the
functor of 3.8 is an equivalence.

{\it Proof.} Let $M$ be a vector space. Suppose we have
an $A_X$-module structure on $i_{x*} M$. We want to show
that it comes from a $T(F_x )\,\tilde{}$-module structure.

(a) The $A^0_X$-action is the same as a chiral
$\ft\tilde{_\CD}$-action = a
$\ft (F_x )\,\tilde{}$-module structure on $M$. Let us show
that our $\ft (F_x )\,\tilde{}$-action is 
integrable, i.e., it comes from a $T(F_x
)^{0}\,\tilde{}$-action. Pick some $\gamma \in\Gamma$; we
want to check that all operators $h_n^\gamma$, $n>0$, act
on $M$ in a locally nilpotent way, and
$h_0^\gamma$ is semi-simple with integral eigenvalues.

(b) We follow the notation from the end of 3.5. For  $\ell$
as in loc.cit.~we have the vertex operator $V^\gamma_\ell
(z):=
\ell\circ\cdot : M
\to M ((z))$ acting on $M$.\footnote{In 3.5 itself we
consider the case $M =A_x$.} 
Then $[h_n^{\gamma'},
V^\gamma_\ell (z)]= c(\gamma ,\gamma' ) z^n V^\gamma_\ell
(z)$ and, by the first equation from (3.6.1), our
$V^\gamma_\ell (z)$ satisfies differential equation
(3.6.2). 

 Consider the modified operator
$\bar{V}^\gamma_\ell (z) :=
\exp (-\mathop\sum\limits_{n>0}h^\gamma_{-n}\frac{z^n}{n})
V^\gamma_\ell (z) :M \to M((z))$. One has:

(i) $\bar{V}^\gamma_\ell$ commutes with every 
$h_{n}^{\gamma'}$, $\gamma'\in\Gamma$, $n>0$.

(ii) $[h_n^{\gamma'},
\bar{V}^\gamma_\ell (z)]=c(\gamma ,\gamma' ) z^n
\bar{V}^\gamma_\ell (z)$ for $n\le 0$.

(iii) $z\partial_z
\bar{V}^\gamma_\ell (z) =
\bar{V}^\gamma_\ell (z)( \mathop\sum\limits_{n\ge 0}
 h_{n}^\gamma z^{-n}) =
(\mathop\sum\limits_{n\ge 0}
 h_{n}^\gamma z^{-n} - c(\gamma ,\gamma
))\bar{V}^\gamma_\ell (z).
$ 

(c) {\it Local nilpotency of
$h^\gamma_{>0}$}: Consider $M$ as a $k[h^\gamma_1,
h^\gamma_2, ..]]$-module. We want to check that it is
supported at
$0\in\Spec k[h^\gamma_1, h^\gamma_2, ..]]:= \cup \Spec
k[h^\gamma_1,..,h^\gamma_n ]$. Take any
$m\in M$. Set $v(z)=\sum v_i z^i :=
\bar{V}^\gamma_\ell (z) m$. Then the support of $v(z)$
(:= the union of supports of $v_i$'s) is equal to that of
$m$. Indeed, Supp $v(z) \subset$ Supp $m$ by (i) above.
To see the opposite inclusion notice that for appropriate
$\ell'
\in\lambda^{-\gamma}$ one has
$ (z-z')^{c(\gamma ,\gamma
)} V^{-\gamma}_{\ell'}(z')V_\ell^\gamma
(z) m=m +
(z-z')\phi (z,z')$ where $\phi \in M[[z',z]][z^{'-1},
z^{-1}]$. Since $V^{-\gamma}_{\ell'}(z')V_\ell^\gamma
(z) m= \exp
(-\mathop\sum\limits_{n>0}h^\gamma_{-n}
\frac{z^{'n}}{n})\bar{V}^{-\gamma}_{\ell'}\exp
(\mathop\sum\limits_{n>0}h^\gamma_{-n}\frac{z^n}{n})v(z)\in
M((z'))((z))$ we see that $m=\sum \alpha_i v_i$ where
$\alpha_i$ are differential operators\footnote{with respect
to the  $k[h^\gamma_1, h^\gamma_2, ..]]$-module structure
on $M$.} acting on $M$ almost all equal to 0. Thus
 Supp $m \subset$ Supp $v(z)$.

Assume that $m$ is killed by $h^\gamma_{>N}$; let
$\eta \in\Spec k[h_1^\gamma ,..,h_N^\gamma ]$ be a generic
point of the support of $m$. We want to show that $\eta
\neq 0$. Changing $m$ if necessary we can assume that the
maximal ideal of $\eta$ kills $m$. Let us localize $M$ (as
a $k[h_1^\gamma ,..,h_N^\gamma ]$-module) at $\eta$. By
above, $v_\eta (z)$ is a non-zero element of $M_\eta
((z))$ killed by the maximal ideal of the local ring.
Differential equation (iii) shows that $z\partial_z
v_\eta (z) = (\mathop\sum\limits_{n\ge 0}
 h_{n\eta}^\gamma z^{-n} - c(\gamma ,\gamma
))v_\eta (z)$ where $h_{n\eta}^\gamma \in k_\eta$ are
images of $h_n^\gamma$ in the residue field. An immediate
induction by $n$ shows that $h_{n\eta}^\gamma
=..=h_{1\eta}^\gamma =0$.

(d) {\it $h^\gamma_0$ is semi-simple with integral
eigenvalues:} By (c) the operator $\tilde{V}^\gamma_\ell
(z)  := \bar{V}^\gamma_\ell
(z)
\exp  (\mathop\sum\limits_{n>0}h_{n}^\gamma
\frac{z^{-n}}{n})$ is well-defined. Since $z\partial_z
\tilde{V}^\gamma_\ell (z) = (h^\gamma_0 -c(\gamma ,\gamma
)) \tilde{V}^\gamma_\ell (z)$ we see that each component 
$\bar{V}^\gamma_n$\footnote{Here $\tilde{V}^\gamma_\ell
(z) =\sum \bar{V}^\gamma_n z^{-n}$.} sends $M$ to the
$c(\gamma ,\gamma )-n$-eigenspace of $h^\gamma_0$. It is
easy to see\footnote{Look at the composition
$V^\gamma_\ell V^{-\gamma}_\ell $ and note that
$\bar{V}^\gamma_\ell$ commutes with $h^\gamma_{\neq 0}$,
cf.~(c).}that the image of
$\bar{V}^\gamma_\ell (z)$ generates $M$, and we are done.

(e) Now we can define the action of $T(F)\,\tilde{}$ on
$M$. We aready know how the connected component
$T(F)^0\,\tilde{}$ acts. We define the action of elements
$\tilde{t}^{-\gamma}$ as $\tilde{V}_\ell^\gamma
(-z)^{-h_0^\gamma}$ (cf.~(3.5.2)); these operators do not
depend on $z$. One check immediately that this provides a
$T(F)\,\tilde{}$-action on
$M$; the corresponding vertex module structure coincides
with the initial one by (3.5.2).   \hfill$\square$

\medskip

{\bf 3.10 Twists and rigidity.} Recall that the category of
chiral algebras on
$X$ is a tensor category. So we know what is a
coalgebra in our tensor category, its (co)action on
some chiral algebra $A_X$, and the chiral subalgebra of
$A_X$ of invariants of this action. We need the case when
our coalgebra is actually (the algebra of functions on) an
affine group $\CD_X$-scheme $\CG_X$. Here we say that
$\CG_X$ acts on $A_X$ and denote the invariants by
$A_X^{\CG_X}$. Suppose now that we have a $\CD_X$-scheme
$\CG_X$-torsor $\CF
=\Spec R_X$ (i.e., a $\CG_X$-torsor  with
connection). The {\it $\CF$-twist} of $A_X$ is, by
definition, $A_X^{\CF}:= (R_X \otimes A_X )^{\CG_X}$. 

\medskip

Our situation: Let
$A$ be a lattice Heisenberg vertex algebra of a {\it
non-degenerate} level
$c$,
$A_X$ the corresponding chiral algebra on $X$. 

(a) Consider the $\Gamma$-grading of $A$ as a
$T^\vee$-action. It yields an action of the constant group
$\CD_X$-scheme $T^\vee_X$ on $A_X$.
So for every de Rham
$T^\vee$-local system $(\fF ,\nabla )$ (:=
a $\CD_X$-scheme $T^\vee_X$-torsor)  we have the twisted
chiral algebra
$A^{(\fF ,\nabla)}_X$ on
$X$.
The $\Gamma$-grading on $A$ yields the one on 
$A^{(\fF ,\nabla)}_X$.

(b) On the other hand,  $T(O)$ acts on $A$ via $i :
T(O)\hra T(F)\,\tilde{}$. The subgroup $Z^c \subset T(O)$
(see 1.9) acts trivially, so, since $c: T(O)/Z \iso
T^\vee (O )$,\footnote{Here we use non-degeneracy of $c$.}
we have a
$T^\vee (O )$-action on
$A$. This provides an action of the group $\CD_X$-scheme
$\CJ_X T^\vee$ of $T^\vee$-valued jets on the chiral
algebra $A_X$. So every $\CD_X$-scheme $\CJ_X
T^\vee$-torsor
$\fF_\CJ$ yields the twisted chiral algebra
$A^{\fF_\CJ}_X$ on $X$. The $T^\vee
(O)$-acton on $A$ preserves the $\Gamma$-grading, so
$A^{\fF_\CJ}_X$ carries a $\Gamma$-grading.

(c) The canonical morphism $\CJ_X T^\vee
\twoheadrightarrow T^\vee_X$ of group $X$-schemes
identifies the sheaf of horizontal sections of $\CJ_X
T^\vee$ with $T^\vee (\CO_X )$. So, by (b), $T^\vee (\CO_X
)$ acts on $A_X$ as a plain sheaf of groups. Thus  any
$T^\vee$-torsor $\fF$ on $X$ yields a $\Gamma$-graded
chiral algebra $A_X^\fF$ := the twist of
$A_X$ by the $T^\vee (\CO_X )$-torsor of sections of $\fF$.

{\it Remarks.} (i) Assume we are in situation (a). Since
$T^\vee$ acts on $A^0$ trivially  we have
a canonical identification of chiral algebras
$$A_X^0 \iso A^{(\fF ,\nabla)0}_X  .\tag 3.10.1$$

(ii) The constructions of (b) and (c) are essentially
equivalent. Indeed, there is a canonical equivalence of
groupoids
$$ \{ T^\vee \text{-torsors on } X \} \iso \{ \CD_X
\text{-scheme }
\CJ_X T^\vee \text{-torsors on } X\} \tag 3.10.2$$
which assigns to a
$T^\vee$-torsor $\fF$ on $X$ the $\CD_X$-scheme $\CJ_X
T^\vee$-torsor of jets $\CJ \fF$; the inverse functor is
the push-out for the canonical homomorphism of the group
$X$-schemes
$\CJ_X T^\vee \twoheadrightarrow T^\vee_X$.  It identifies
the $T^\vee (\CO_X )$-torsor of sections of $\fF$ with that
of horizontal sections of
$\CJ\fF$, so
 one has a canonical identification 
$$ A^{\fF}_X \iso  A^{\CJ\fF}_X . \tag 3.10.3$$

{\bf 3.11 Lemma.} (i) For every de Rham $T^\vee$-local
system $(\fF ,\nabla )$ there is a canonical isomorphism of
chiral algebras 
$$ \phi :A^{(\fF ,\nabla )}_X \iso A^\fF_X . \tag 3.11.1$$ 

(ii) A section  $s\in\fF$ yields $\nu := \nabla\log
s \in \ft^\vee \otimes \omega (X)$ and an
identification of chiral algebras $\beta_s : A_X
\iso A^\fF_X$. Now (3.10.1) equals the composition of
$\phi^{-1}\beta_s$  with the automorphism of $A^0_X =
U(\ft_\CD )$ coming from an automorphism 
 $\tilde{a} \mapsto \tilde{a}+ (\nu ,a)$ of
$\ft\tilde{_\CD}$.

{\it Proof.} Key remark: the action of $T^\vee$ from
3.10 (a) coincides, via  $T^\vee
\subset T^\vee (O)$, with the action from 3.10 (b).
Therefore $A^{(\fF ,\nabla )}_X$ coincides with the twist
of $A_X$ by the $\CD_X$-scheme $\CJ_X T^\vee$-torsor
induced from $(\fF ,\nabla )$ by the embedding of group
$\CD_X$-schemes $T^\vee_X \hra \CJ_X T^\vee$. By (3.10.2)
the latter torsor equals $\CJ \fF$, so (i) comes
from (3.10.3). 

In more details, consider the
canonical morphism of group $X$-schemes $\CJ_X T^\vee
\twoheadrightarrow T^\vee_X$ and the one of torsors $\CJ\fF
\twoheadrightarrow\fF$. The ``constant" connection on
$T^\vee_X$ yields a section $T^\vee_X \hra \CJ_X T^\vee$
which is an embedding of group $\CD_X$-schemes, and 
$\nabla$ yields a section $e_\nabla :\fF \hra \CJ\fF$ which
is an embedding of $\CD_X$-schemes compatible with the
actions of $T^\vee_X \hra \CJ_X T^\vee$. This $e_\nabla$
identifies
$\CJ\fF$ with the $\CD_X$-scheme $\CJ_X T^\vee$-torsor
induced from $(\fF ,\nabla )$ via the embedding of the
group $\CD_X$-schemes $T^\vee_X \hra \CJ_X T^\vee$. Our
$\phi$ is the composition 
$ A^{(\fF ,\nabla )}_X
\buildrel{e_\nabla}\over\longrightarrow A^{\CJ\fF }_X
\buildrel{(3.10.3)}\over\longleftarrow A^\fF_X . $ 

Let us write down $\phi$ in terms of  $s$ as in (ii). Let
$\CJ s$ be the horizontal section of $\CJ\fF$ defined by
$s$. The projection $\CJ\fF \twoheadrightarrow \fF$ maps
$\CJ s$ to $s$, hence $$e_\nabla (s)= \kappa_\nu \CJ s
\tag 3.11.2$$ where
$\kappa_\nu $ is a section of $\CJ_X
T^\vee$ killed by $\CJ_X T^\vee
\twoheadrightarrow T^\vee_X$ uniquely determined by the
condition
$d\log
\kappa_\nu =\nu \in \ft^\vee
\otimes\omega_X \subset \CJ_X \ft^\vee \otimes \omega_X$. 

As a mere
$\CO_X$-module $A_X^{(\fF,\nabla )}$ is the
$\fF$-twist of
$A_X$ with respect to the action of $T^\vee$ given by the
$\Gamma$-grading (see 3.10 (a)), so
$s$ yields an identification  $\alpha_s :A_X
\iso A^{(\fF ,\nabla )}_X$ of $\CO_X$-modules. Similarly,
$\CJ s$ yields an identification of chiral
algebras $\alpha_{\CJ s} : A_X
\iso A^{\CJ\fF}_X
$ which equals the composition of $\beta_s$ (see
3.11(ii)) and (3.10.3). By (3.11.2) one has
$$\phi
\alpha_s =
\beta_s \kappa_\nu
 : A_X \iso A_X^\fF . \tag 3.11.3$$

On $A^0_X =U(\ft\tilde{_\CD} )$ our $\alpha_s$ equals
(3.10.1) and
$\kappa_\nu$ acts  according
to the automorphism $\tilde{a} \mapsto \tilde{a} + (\nu ,a
)$ of $\ft\tilde{_\CD}$ (see 1.5(ii)). This implies (ii).
\hfill$\square$

\bigskip

{\bf \S4 Moving $T^\vee$-local system}  

\medskip
 
{\bf 4.1} Let 
$\CL\CS =\CL\CS_{T^\vee}$ be the moduli stack of de
Rham $T^\vee$-local systems on Spec$\, F$. By definition,
an $R$-point of $\CL\CS$ is a pair $(\fF_R
,\nabla )$ where $\fF_R$ is a $T^\vee$-torsor on
$\Spec F_R$, $F_R := F\hat{\otimes} R =\limleft
F/\fm^n
\otimes R \simeq R((t))$, and $\nabla$ is an $R$-relative
continuous connection on $\fF_R$. We always assume that 
$\fF_R$ comes from a $T^\vee$-torsor on $\Spec O_R$, $O_R
:= O\hat{\otimes}R \simeq R[[t]]$, \'etale locally on
$\Spec R$.  

According to Drinfeld, the latter unpleasant
asssumption becomes redundant if we consider instead of
\'etale topology a finer one (a version of cdh topology),
but this seems to be irrelevant for what follows.

{\bf 4.2} Here are
some convenient descriptions of $\CL\CS$: 

 Set $\omega (F)
:=\{$1-forms on Spec$F\}$; this is an ind-scheme.
The ind-scheme $\CC$ of connections on the trivialized
$ T^\vee$-bundle  identifies canonically with $ 
\ft^\vee \otimes\omega (F)$,$\nu\mapsto
\nabla_\nu := \partial_t +\nu$. Notice that $$\CC = \Spec
(\Sym\, \ft (F) ):= \limright \Spec (\Sym\, \ft (F)/\ft
(\fm^n ))
\tag 4.2.1$$ where $a\in \ft (F)$ is identified with a
linear function $\phi (a): \nu \mapsto \Res (a,\nu )$ on
$\CC$.

The group ind-scheme $
T^\vee (F)$ of automorphisms of a
$ T^\vee$-bundle acts on $\CC$;
the corresponding gauge action on $\ft^\vee \otimes\omega
(F)$ is
$g(\nu )= \nu +d\log ( g )$. Thus
$$\CL\CS = \CC/T^\vee (F)=\ft^\vee \otimes\omega (F)
/T^\vee (F).
\tag 4.2.2$$ 

 Set $\Phi :=T^\vee (F)
/T^\vee (\fm )$. Thus $T^\vee =T^\vee (O)/T^\vee (\fm )$ is
a subgroup of
$
\Phi$ and 
$\Phi / T^\vee = Q\times  \Gamma^\vee $ where $Q$
is the formal group whose Lie algebra equals $ \ft^\vee
\otimes (F/O)$. Set
$\omega (F)^{-}:=
\omega (F) /\omega (O)$, $\bar{\CC} := \ft^\vee \otimes
\omega (F)^- $; one has $\bar{\CC}=\Spec (\Sym \,\ft (O)
):=\limright \Spec (\Sym \,\ft (O/\fm^n ))$ (see (4.2.1)).
Since
$d\log$ yields an isomorphism
$ T^\vee (\fm )
\iso
\ft^\vee \otimes \omega (O) $ one has $\bar{\CC}=\CC
/T^\vee (\fm )$ and
$$\CL\CS = 
 \bar{\CC} /\Phi . \tag 4.2.3$$ 

There is a canonical decomposition $\omega (F)^- =
\omega (F)^{irr}\times k$ where
the projection $
\omega (F)^{-}\to k$ is the residue map, and $k\hra
\omega (F)^{-}$ is the subspace of forms with pole of
order one. Hence $\bar{\CC} =\bar{\CC}^{irr}
\times\ft^\vee$. The morphism
$d\log :\Phi\to
 \ft^\vee \otimes \omega (F)^- $ kills
$ T^\vee$, so it yields a morphism $ Q\times \Gamma \to
(\ft^\vee
\otimes
\omega (F)^{irr})\times  \ft^\vee$. This is an embedding
compatible with the product decomposition; it identifies
$Q$ with the formal completion of $\ft^\vee \otimes
\omega (F)^{irr}$, and $\Gamma^\vee \hra  \ft^\vee$ is
the usual embedding. Therefore, by (4.2.3), we have a
canonical projection $$\CL\CS \to (\bar{\CC}^{irr}/Q)\times
(\ft^\vee /\Gamma^\vee )  \tag 4.2.4$$ which makes
$\CL\CS$ a
$T^\vee$-gerbe over $(\bar{\CC}^{irr}/Q)\times (t^\vee
/\Gamma^\vee )$. Any splitting  of the extension $0\to
T^\vee \to\Phi
\to Q\times \Gamma \to 0$ yields a trivialization of this
gerbe, i.e., an identification $$ \CL\CS
\iso (\bar{\CC}^{irr}/Q)\times (t^\vee /\Gamma^\vee )
\times B T^\vee .
\tag 4.2.5$$ 

{\bf 4.3} The definition of 3.10(a) works for families of
de Rham $T^\vee$-local systems, so we have an
$\CL\CS$-family $A^{\CL\CS}_{F}$ of twisted chiral
algebras  over
$\Spec F$ equipped with a  $\Gamma$-grading. 

Namely, recall that a $\Spec R$-family of chiral
algebras on $\Spec F$ is a flat $F_R$-module $A_{F_R}$
equipped with an
$R$-relative connection and a chiral product operation
which amounts to an $R$-bilinear ope product $
A_{F_R} \otimes A_{F_R} \to A_{F_R \, 1 }((t_1 -
t_2 ))$ which satisfies the usual properties.

For $(\fF_R ,\nabla )\in \CL\CS (R)$ one defines the
$(\fF_R ,\nabla )$-twisted chiral lattice algebra
$A_{F_R}^{(\fF_R ,\nabla )}$ as in 3.10(a). This
construction is compatible with the $R$-base change. The
 datum of $A_{F_R}^{(\fF_R ,\nabla )}$ together with
the base change isomorphisms forms
$A^{\CL\CS}_{F}$.

The pull-back of $A^{\CL\CS}_{ F}$ to $\CC
=\ft^\vee \otimes \omega (F)$ is a $\CC$-family of chiral
algebras $A^{\CC}_{ F}$ on
$\Spec F$ equivariant with respect to the action of the
group ind-scheme $T^\vee (F)$.

{\bf 4.4} Denote by $\CM (A^{\CL\CS}_{ F} )$ the
abelian
$k$-category of $A^{\CL\CS}_{F}$-modules, i.e., of
$\CO^!$-modules on $\CL\CS$ equipped with an action of
$A^{\CL\CS}_{F}$. Explicitly, an $A^{\CL\CS}_{F}$-module
$M$ amounts to a disctrete module
$M_\CC$ over the topological algebra $\CO (\CC )=\CO
(\ft^\vee
\otimes \omega (F))=\limleft \Sym_k (\ft \otimes F/\fm^n )$
equipped with compatible actions of the group
ind-scheme
$T^\vee (F)$ and  the chiral algebra $A^{\CC}_{F}$.

The vector space of sections $\Gamma (M):=\Gamma (\CL\CS
,M)$ is defined as a BRST-type combination of
invariants and coinvariants: $$\Gamma (M)
:= (M_\CC^{T^\vee (O_F )})_{T^\vee (F)/T^\vee (O_F )}. \tag
4.4.1$$

Since $A_F^{\CL\CS 0}\subset A_F^{\CL\CS}$ is a
$\CC$-constant family of chiral algebras equal
$U(\ft\tilde{_\CD})_\CC$ (see (3.10.1)) the Heisenberg Lie
algebra $\ft (F)\,\tilde{}$ acts on $M_\CC$ and
$\Gamma (M)$.

\medskip

{\bf 4.5 Theorem.} The functor $\Gamma :\CM (A^{\CL\CS}_{F}
)\to \ft (F)\,\tilde{}$-mod is an equivalence of
categories.

{\it Proof.} (a) Recall that $T^\vee (F)$ acts on
$T(F)\,\tilde{}$ (see 1.5(iii)); denote by $(T^\vee \times
T)(F)\,\tilde{}$ the corresponding semi-direct product.
This is a symmetric
Heisenberg extension for the torus $T^\vee \times T$ whose
level is of  determinant 1.  According to 1.9 the
functor $ P \mapsto  P^{T^\vee (O)\times
T(O)}$ from the category of $(T^\vee \times
T)(F)\,\tilde{}$-modules to (super)
vector spaces is an equivalence of categories.\footnote{The
inverse functor is induction from $T^\vee (O)\times T(O)$
to $(T^\vee (F)\times T(F))\,\tilde{}$.}

\medskip

(b) The group $(T^\vee \times T)(F)\,\tilde{}$ acts on
$\CC$ via the projection to $T^\vee (F)$, so we have the
category $\CP$ of $(T^\vee \times
T)(F)\,\tilde{}$-equivariant discrete $\CO (\CC
)$-modules.\footnote{As always,  the center
$\Bbb G_m$ of $(T^\vee \times
T)(F)\,\tilde{}$ acts on our modules by standard
homotheties.}

{\bf Lemma.} (i) For every $M\in \CM (A^{\CL\CS}_{F})$
there is a natural $T(F)\,\tilde{}$-action on $M_\CC$
which, together with the structure $\CO (\CC )$- and
$T^\vee (F)$-actions, makes $M_\CC$ an object of $\CP$. 
The functor $\CM (A^{\CL\CS}_{F}) \to\CP$ is an
equivalence of categories.

(ii) The $\ft (F)\,\tilde{}$-action on $M_\CC$ from 4.4 
is $\tilde{a}m = \alpha (\tilde{a})\star m$
where $\star$ is the $\ft (F)\,\tilde{}$-action coming from
the $T(F)\,\tilde{}$-module structure from (i) and
$\alpha $ is a $\CC$-family of automorphisms of $\ft
(F)\,\tilde{}$, $\alpha (\tilde{a}) := \tilde{a} +
\phi (a)$ (see (4.2.1)).

{\it Proof of Lemma.} 
Notice that Lemma 3.11, as well as its proof,
remains valid for families of de Rham $T^\vee$-local
systems; one can also replace the curve $X$ by
$\Spec F$. By 3.11(i) our $\CC$-family of
chiral algebras 
$A^{\CC}_{F} $ on $\Spec F$ can be canonically identified
 with the constant $\CC$-family
$A_{F \CC}$. This isomorphism is compatible with the
$T^\vee (F)$-actions (here the $T^\vee (F)$-action on
$A_{F \CC}$ is the product of the action on $A_F$ from
3.10(c) and the gauge action on $\CC$).  

Therefore an $A^{\CL\CS}_F$-module $M$ amounts to a
 $T^\vee (F)$-equivariant discrete $\CO (\CC )$-module
$M_\CC$ equipped with an
$A_F$-action which commutes with the $\CO (\CC )$-action
and is compatible with the $T^\vee (F)$ one. According to
3.8, 3.9 an $A_F$-action is the same as a
$T(F)\,\tilde{}$-action. As follows from the construction
of 3.8, the compatibility with the $T^\vee (F)$-action just
means that the $ T^\vee (F)$- and $
T(F)\,\tilde{}$-actions form an action of $(T^\vee \times
T)(F)\,\tilde{}$. This proves (i), and (ii) follows
from 3.11(ii).   \hfill$\square$

\medskip

(c) Set $H:= T(F)^0 /T$, $H^\vee := T^\vee (F)^0 /T^\vee$.
Since $T$ is a central subgroup of $T(F)^0\,\tilde{}$ and
$T^\vee \subset T^\vee (F)^0$ acts trivially on
$T(F)^0\,\tilde{}$ we have a central $\Bbb G_m$-extension
$H\,\tilde{}:= T(F)^0\,\tilde{}/T$ of $H$ and the
corresponding semi-direct product $(H^\vee \times
H)\,\tilde{}$ of $H^\vee $ and $H\,\tilde{}$.
Since $T^\vee (F)^0$ acts on $\CC$ through its quotient
$H^\vee$ we have the
category $\CP^0$ of $ (H^\vee \times H)
\,\tilde{}$-equivariant discrete $\CO (\CC
)$-modules. As follows  from (a)  the
obvious functor $\CP \to \CP^0$, $M_\CC \to M_\CC^0 :=
(M_\CC )^{T^\vee \times T}$, is an equivalence of
categories.

\medskip

(d) As we have seen in 4.2, $\CC$ is a 
$T^\vee (\fm )$-torsor over $\bar{\CC}$
(see 4.2). The
$T^\vee (\fm )$-action on $H \,\tilde{}$
defines the corresponding twisted group ind-scheme
 on $\bar{\CC}$ which we denote by
$H\tilde{_{\bar{\CC}}}$; this is a central $\Bbb
G_m$-extension of the constant group ind-scheme
$H_{\bar{\CC}}$. The action of
$H^\vee $ on $\CC$ and
$H\,\tilde{}$ defines an action of the formal group
ind-scheme $Q:=T^\vee (F)^0 /T^\vee (O) =H^\vee /T^\vee
(\fm )$ on
$\bar{\CC}$ and
$H\tilde{_{\bar{\CC}}}$.

Let $\bar{\CP}$ be the category of discrete $\CO
(\bar{\CC})$-modules $M_{\bar{\CC}}$ equipped with an
$H\tilde{_{\bar{\CC}}}$-action and equivariant with
respect to the $Q$-action. There is an obvious descent
equivalence  $\CP^0 \iso \bar{\CP}$, $N_{\CC}
\mapsto N_{\bar{\CC}}:= N_\CC^{T^\vee (\fm )}$.

\medskip

(e) The $\CC$-family $\alpha$ of automorphisms of $\ft
(F)\,\tilde{}$ from (b)(ii) has property $g (\alpha (
\tilde{a})_\nu )= \alpha (\tilde{a}_{g\nu} )$ for every
$g\in T^\vee (F)$, $\nu \in \CC $, $\tilde{a}\in\ft
(F)\,\tilde{}$. Therefore it yields a canonical
morphism of Lie algebras $\bar{\alpha} : \ft
(F)\tilde{_{\bar{\CC}}} \to
\fh\tilde{_{\bar{\CC}}}$ on
$\bar{\CC}$ equivariant with respect to the $Q$-action.
Here $\ft (F)\tilde{_{\bar{\CC}}}$ is the
constant $\bar{\CC}$-family and $\fh\tilde{_{\bar{\CC}}}
$ is the Lie algebra of $H\tilde{_{\bar{\CC}}} $. Thus
every $N_{\bar{\CC}}\in\bar{\CP}$ is a $\ft
(F)\,\tilde{}$-module in a natural way. The $\ft
(F)\,\tilde{}$-action commutes with the $Q$-action, so
the space of $Q$-coinvariants $(N_{\bar{\CC}})_Q$ is also
a $\ft
(F)\,\tilde{}$-module.

We have an equivalence $\CM (A^{\CL\CS}_F )\iso
\bar{\CP}$, $M\mapsto M^0_{\bar{\CC}}$, defined as the
composition of the equivalences from (b)(i), (c), and (d).
By (4.4.1), (b)(ii) there is a canonical identification of $\ft
(F)\,\tilde{}$-modules
$\Gamma (M)= (M^0_{\bar{\CC}})_Q$. 
Therefore our theorem amounts to the fact that the functor
$\bar{\CP} \to \ft (F)\,\tilde{}$-mod,
$N_{\bar{\CC}}\mapsto (N_{\bar{\CC}})_Q$, is an
equivalence of categories.

\medskip

(f) The morphism $\bar{\alpha} :\ft
(F)\tilde{_{\bar{\CC}}} \to
\fh\tilde{_{\bar{\CC}}}$ is surjective and
$H\tilde{_{\bar{\CC}}}$ is connected, so for 
$N_{\bar{\CC}}\in\bar{\CP}$ the
$H\tilde{_{\bar{\CC}}}$-action on it is uniquely
determined by the $\ft (F)\,\tilde{}$-action. Therefore
$\bar{\CP}$ is a full subcategory of the category  of 
$Q$-equivariant discrete $\CO (\bar{\CC})$-modules
equipped with a $\ft (F)\,\tilde{}$-action. 

Notice that $H$ is an extension of a formal group by a
unipotent one $T(\fm )$. Thus $H\,\tilde{}$-modules are
the same as $\fh\,\tilde{}$-modules such that the Lie
algebra of $T(\fm )\subset H\,\tilde{}$
acts in a locally nilpotent way. Same is true for
$H\tilde{_{\bar{\CC}}}$-modules. 

Let $\beta : \ft (O)_{\bar{\CC}}\to \ft
(F)\tilde{_{\bar{\CC}}}$ be the morphism $\beta (a):=
a-\phi (a)$.\footnote{Here $\phi (a)\in \CO (\bar{\CC})$
is the linear function on $\bar{\CC}$ that corresponds to
$a$, see 4.2.} Then $\Ker \bar{\alpha} =\beta
(\ft_{\bar{\CC}})$ and $\bar{\alpha}$ identifies 
$\beta (\ft (\fm )_{\bar{\CC}})$ with the Lie algebra
of $T (\fm )\subset H\tilde{_{\bar{\CC}}}$.

We see that objects of $\bar{\CP}$ are exactly those
$Q$-equivariant discrete $\CO (\bar{\CC})$-modules
$N_{\bar{\CC}}$ equipped with a $\ft (F)\,\tilde{}$-action
which satisfy the following conditions: 

(i) $\beta (\ft )$ kills
$N_{\bar{\CC}}$,

(ii) $\beta (\ft (\fm )_{\bar{\CC}})$ acts on
$N_{\bar{\CC}}$ in a locally nilpotent way.

\medskip

(g) It is convenient to view $\ft
(F)\,\tilde{}$-modules in the following ``geometric" way.

According to 4.2 the embedding $\phi : \ft
(O)\hra \CO (\bar{\CC})$ identifies $\bar{\CC}$ with $\Spec
(\Sym \ft (O))$. Let $\pi :\CC
\to \ft^\vee$ be the projection that corresponds to the
embedding $\ft\hra
\ft (O)$.\footnote{ $\pi$ can be also
described as the residue map $\CC =\ft^\vee \otimes
(\omega (F)/\omega (O)) \to \ft^\vee$, see 4.2.} The Lie
algebra
$\ft (F)$ acts along the fibers of $\pi$: namely,
for $a \in \ft (F)$, $b\in \ft (O)\subset \CO (\CC
)$ one has $a(b):= \Res c(b,da)$. Let
$\Theta\tilde{_{\pi}}$ be a Lie algebroid on $\bar{\CC}$
equipped with a Lie algebra homomorphism $\rho : \ft
(F)\,\tilde{} \to \Theta\tilde{_{\pi}}$ such that:

(i) $\Theta\tilde{_{\pi}}$ is an extension of the relative
tangent algebroid $\Theta_\pi
=\Theta_{\bar{\CC}/\ft^\vee}$ by $\CO_{\bar{\CC}}$,

(ii) $\rho$ lifts the above action of $\ft (F)$ on
$\bar{\CC}$, and $\rho |_{\ft (O)}$ equals $\phi
:\ft (O)\to\CO (\bar{\CC})\subset \Theta\tilde{_{\pi}}$.

Such $\Theta\tilde{_\pi}$ exists and is unique.

By definition, a $\Theta\tilde{_\pi}$-module any discrete
$\CO (\bar{\CC})$-module $L_{\bar{\CC}}$ equipped with a
(right)\footnote{Since $\bar{\CC}$ is an ind-scheme only
right actions make sense.} action of
$\Theta\tilde{_\pi}$ such that
$1\in \CO_{\bar{\CC}}\subset
\Theta\tilde{_\pi}$ acts as identity. Then
 $L_{\bar{\CC}}$ is a $\ft (F)\,\tilde{}$-module via
$\rho$. The functor $\Theta\tilde{_\pi}$-mod $\to\ft
(F)\,\tilde{}$-mod is an equivalence of categories.

So  $\ft (F)\,\tilde{}$-modules are the same as (twisted)
$\CD$-modules along the fibers of $\pi$.

\medskip

(h) Now we are ready to prove the theorem. Returning to
(f) notice that $Q$ is a formal group which acts formally
simply transitively along the fibers of $\pi$. Therefore a
$Q$-equivariant structure on an $\CO (\bar{\CC} )$-module
is the same as a (right) $\Theta_\pi$-action. Such object
equipped with an $\ft (F)\,\tilde{}$-action amounts, by
(g), to a $\Theta_\pi \times
\Theta\tilde{_\pi}$-module, i.e., a discrete $\CO
(\bar{\CC}\times\bar{\CC})$-module
$N_{\bar{\CC}\times\bar{\CC}}$ equipped with a (right)
action of $\Theta_\pi \times
\Theta\tilde{_\pi}$ such that $1 \in \CO_{\bar{\CC}}\subset
\Theta\tilde{_\pi}$ acts as identity. Condition (f)(i)
means that $N_{\bar{\CC}\times\bar{\CC}}$ is supported
(scheme-theoretically) on
$\bar{\CC}\mathop\times\limits_{\ft^\vee}\bar{\CC} \subset 
\bar{\CC}\times\bar{\CC}$. Condition (f)(ii) means that
$N_{\bar{\CC}\times\bar{\CC}}$ is supported
set-theoretically on the diagonal $\bar{\CC}\subset
\bar{\CC}\times\bar{\CC}$.

According to the Kashiwara lemma such animals are the same
as $\Theta\tilde{_\pi}$-modules on the diagonal $\bar{\CC}\subset
\bar{\CC}\times\bar{\CC}$. The identification can be given
by the integration functor along the first copy of
$\bar{\CC}$ which is the same as taking $Q$-coinvariants. 
Returning from $\Theta\tilde{_\pi}$-modules to $\ft
(F)\,\tilde{}$-modules we get the theorem.
\hfill$\square$

 {\bf References.}

\medskip

[BBE] A.~Beilinson, S.~Bloch, H.~Esnault. $\CE$-factors
for Gau\ss-Manin determinants. alg-geom 2001.

[BD] A.~Beilinson, V.~Drinfeld. Chiral Algebras.

[CC] C.~E.~Contou-Carr\`ere.
Jacobienne locale, groupe de
bivecteurs de Witt universel et symbole local
mod\'er\'e. C.~R.~Acad.~Sci.~Paris,
t.~318, S\'erie I (1994), 743-746.

[D] C.~Dong. Vertex algebras  associated with even
lattices. J.~Algebra 161 (1993) 245--265.

[FBZ] E.~Frenkel, D.~Ben-Zvi. Vertex algebras and
algebraic curves. AMS 2001.

[K] V.~Kac. Vertex algebras for beginners. University
Lecture Series vol.~10, AMS, 1998.

\end